\documentclass{amsart}
\pdfoutput=1

\usepackage{amssymb}
\usepackage{microtype}
\usepackage{tikz}
\usetikzlibrary{arrows}
\usepackage{booktabs}
\usepackage[pdftex,colorlinks,citecolor=black,linkcolor=black,urlcolor=black,bookmarks=false]{hyperref}

\usepackage{doi}
\usepackage{subfigure}
\usepackage{bbm} 

\newtheorem{theorem}{Theorem}[section]
\newtheorem{proposition}[theorem]{Proposition}
\newtheorem{lemma}[theorem]{Lemma}
\newtheorem{corollary}[theorem]{Corollary}

\theoremstyle{definition}

\numberwithin{figure}{section}
\numberwithin{equation}{section}
\numberwithin{table}{section}

\newcommand{\Z}{\mathbb{Z}}

\newcommand{\C}{\mathbb{C}}
\newcommand{\F}{\mathbb{F}}

\newcommand{\GL}{{\rm GL}}
\newcommand{\PGL}{{\rm PGL}}
\newcommand{\SL}{{\rm SL}}

\newcommand{\Stab}{{\rm Stab}}

\newcommand{\Orb}{{\rm Orb}}

\DeclareMathOperator{\tr}{tr}
\newcommand{\legendre}[2]{\genfrac{(}{)}{}{}{#1}{#2}}

\title{Experiments with the Markoff Surface}

\author{Matthew de Courcy-Ireland}
\address{Department of Mathematics\\
Princeton University\\
Princeton NJ 08544} \email{mdc4@math.princeton.edu}

\author{Seungjae Lee}
\address{Department of Mathematics\\
Princeton University\\
Princeton NJ 08544} \email{seungjl@math.princeton.edu}

\date{December 18, 2018}

\begin{document}

\begin{abstract}
We confirm, for the primes up to 3000, the conjecture of Bourgain, Gamburd, and Sarnak on strong approximation for the Markoff surface $x^{2}+y^{2}+z^{2} = 3xyz$
modulo primes. 
For primes congruent to 3 modulo 4, we find data suggesting that some natural graphs constructed from this equation are asymptotically Ramanujan. For primes congruent to 1 modulo 4, the data suggests a weaker spectral gap. In both cases, there is close agreement with the Kesten-McKay law for the density of states for random 3-regular graphs. We also study the connectedness of other level sets $x^{2}+y^{2}+z^{2}-3xyz = k$.
In the degenerate case of the Cayley cubic, we give a complete description of the orbits.
\end{abstract}

\maketitle

\section{Introduction} \label{sec:introduction}

The \emph{Markoff equation} or \emph{Markoff surface}
\begin{equation} \label{eqn:markoff}
x^2+y^2+z^2=3xyz
\end{equation}
is preserved by the operations 
\begin{equation}
m_1: x \mapsto 3yz-x, m_2: y \mapsto 3xz-y, m_3: z \mapsto 3xy-z.
\end{equation}
Indeed, (\ref{eqn:markoff}) is a cubic equation overall, but only quadratic in each individual variable, and these moves amount to switching to the other root of the quadratic. They are called \emph{Markoff moves} or \emph{Vieta operations}.
Markoff proved in 1880 \cite{M} that any solution to (\ref{eqn:markoff}) in nonnegative integers except $(0,0,0)$ can be reached by a sequence of Vieta operations and transpositions. This means that all solutions to equation~(\ref{eqn:markoff}) can be found quickly by navigating from the root $(1,1,1)$ using the moves $m_1,m_2,m_3$. This can be represented graphically as a tree with a vertex for each solution, and edges between solutions that are connected by one of the moves $m_1,m_2,m_3$.

For the Markoff equation over a prime field $\F_p$, it is no longer guaranteed that all solutions can be found by these moves. Does every solution mod $p$ lift to a solution over the integers? If so, then the same sequence of Vieta moves used to reach the lift will reach its image mod $p$ because the moves are polynomial operations in $(x,y,z)$.  Over a finite field, instead of the Markoff tree we have a \emph{Markoff graph} with cycles. 
To streamline the graphs, we found it convenient to use the following operations known as \emph{Dehn twists} instead of the Markoff moves $m_1, m_2, m_3$.
The three operations $D_1,D_2,D_3$ are
\begin{align} \label{eqn:dehn}
D_1 : (x,y,z) \mapsto (3yz-x, z, y) \\ 
D_2 : (x,y,z) \mapsto (x, z, 3xz-y) \\ 
D_3 : (x,y,z) \mapsto (x, 3xy-z, y)
\end{align}
With either choice of generators, one never has parallel edges because $(0,0,0)$ is the only solution to $m_j m_k(x,y,z) = (x,y,z)$, and likewise the only solution to $D_jD_k^{-1} (x,y,z) = (x,y,z)$ for $j \neq k$.
Whereas each Markoff move has on the order of $p$ fixed points, the Dehn twists have only a bounded number. Indeed, $D_2$ and $D_3$ have no fixed points except $(0,0,0)$ solving (\ref{eqn:markoff}). If $D_1(x,y,z) = (x,y,z)$, then $z = y$ and $2x = 3y^2$. Substituting into Markoff's equation gives $9y^4/4+2y^2=9y^4/2$. So the only solutions are $(0,0,0)$ and, if $8$ is a square mod $p$, $(4,\pm \sqrt{8},\pm \sqrt{8})/3$.
For each point $(x,y,z) \neq (0,0,0)$ on the Markoff surface over $\F_p$, we take an edge between $(x,y,z)$ and each of its images under $D_1, D_2, D_3$. This defines a 3-regular graph with at most two loops for each prime $p$, often no loops at all, and in any case no parallel edges. Note that $D_j = \tau_{23} \circ m_j$, where $\tau_{23}$ is the transposition exchanging the second and third coordinates.

Note that if $x^2+y^2+z^2=3xyz$, then $(X,Y,Z) = (x,y,z)/3$ solves 
\begin{equation} \label{eqn:xyz}
X^2 + Y^2 + Z^2 = XYZ.
\end{equation}
Over the integers, the factor $3xyz$ in (\ref{eqn:markoff}) guarantees that the base solution is $(1,1,1)$ instead of $(3,3,3)$.
Over $\F_p$ with $p > 3$, it can be more convenient to use version (\ref{eqn:xyz}) of the Markoff surface. 
The corresponding Markoff moves are $x \mapsto yz-x, y \mapsto xz-y, z \mapsto xy-z$. We will denote these also by $m_1, m_2, m_3$ as it will be clear from context whether we are using (\ref{eqn:markoff}) or (\ref{eqn:xyz}).

The connectedness of these graphs for all $p$ is the question of whether \emph{strong approximation} holds for the Markoff surface.
Baragar was the first to conjecture that this connectedness does hold for all $p$ and verified it for $p \leq 179$ (see p. 124 of \cite{Bar}). The present paper extends this to $p < 3000$ and suggests that the graphs are not only connected but even form an expander family as $p$ grows. 
In Section~\ref{sec:zerok}, we present numerical evidence that the Markoff graphs have a spectral gap. For $p \equiv 3 \bmod 4$, the gap is almost as large as possible, while for $p \equiv 1 \bmod 4$ it is somewhat weaker (Figure~\ref{fig:eigenvalue}). We also discuss the bulk distribution of eigenvalues. The data suggest that this converges to the Kesten-McKay law (Figure \ref{fig:histogram}), which has recently been confirmed theoretically by Magee and de Courcy-Ireland \cite{DM}.

In addition to the numerical evidence, there are compelling theoretical reasons to believe in strong approximation. Bourgain-Gamburd-Sarnak proved that it holds unless $p^2-1$ has many prime factors, which happens only for rare values of $p$ \cite{BGS}. 
The strong approximation conjecture is equivalent to a certain group action being transitive, namely the action of Vieta moves and coordinate permutations on solutions $(x,yz)$ modulo $p$ to equation~(\ref{eqn:markoff}), up to double sign changes $(x,y,z) \mapsto (\sigma_1 x, \sigma_2 y, \sigma_3 z)$ with each $\sigma_j = \pm 1$ and $\sigma_1 \sigma_2 \sigma_3 = 1$. Meiri and Puder show that, if $p \equiv 1$ mod $4$ and $p^2 -1$ is not very smooth, so that the analysis of Bourgain-Gamburd-Sarnak shows the action to be transitive, then the resulting permutation group is either the alternating group or the symmetric group on the set of blocks (modulo sign change) \cite{MPC}. 
They prove this for $p \equiv 3$ mod 4 as well, but this requires an additional hypothesis about $p$. Carmon shows in the appendix to \cite{MPC} that this assumption holds except for another sparse sequence of primes. It had been conjectured around the same time by Cerbu-Gunther-Magee-Peilen that the group is alternating for $p \equiv 3$ mod 16 and symmetric otherwise \cite{CGMP}.   
This phenomenon of being fully transitive is a kindred spirit to expansion: Not only is the graph/action connected/transitive, but it is very robustly connected so that there are many ways to go from one point to any other.

\begin{figure}
\centering
\begin{subfigure}{}
\begin{tikzpicture}[scale=1.5]
\draw (0,0) -- (0,1);
\draw (0,0) -- (-0.8660254, -1/2);
\draw (0,0) -- (0.8660254, -1/2);
\draw (0.325,0.2) node {333};
\draw (0.325,0.825) node {633};
\draw (1.325, 1.5) node {(6,15,3)};
\draw (0.625, 2) node {(39,15,3)};
\draw (1.75, 1) node {(6,15,87)};
\draw (0,1) -- (-3*0.8660254/4,1+2/3*1/2);
\draw (0,1) -- (3*0.8660254/4, 1+2/3*1/2);
\draw (-0.8660254,-1/2) -- (-0.8660254,-1/2-2/3);
\draw (-0.8660254,-1/2) -- (-0.8660254 - 3*0.8660254/4, -1/2+2/3*1/2);
\draw (0.8660254,-1/2) -- (0.8660254,-1/2-2/3);
\draw (0.8660254,-1/2) -- (0.8660254 + 3*0.8660254/4, -1/2+2/3*1/2);
\draw (3*0.8660254/4,1+2/3*1/2) -- (3*0.8660254/4,1+2/3*1/2 +3*2/3/4) ;
\draw (3*0.8660254/4,1+2/3*1/2) -- (3*0.8660254/4 + 2/3*2/3*0.8660254 ,1+2/3*1/2 - 2/3*2/3*1/2) ;
\end{tikzpicture}
\end{subfigure}
\\
\begin{subfigure}{}
\begin{tikzpicture}[scale=0.50]
\draw (0,0) -- (0,1);
\draw (0,0) -- (-0.8660254, -1/2);
\draw (0,0) -- (0.8660254, -1/2);
\draw (0,1) -- (-0.8660254, 3/2);
\draw (0,1) -- (0.8660254, 3/2);
\draw (0.8660254, 3/2) -- (0.8660254, 5/2);
\draw (-0.8660254, 3/2) -- (-0.8660254, 5/2);
\draw (0.8660254, 5/2) -- (0,3);
\draw (-0.8660254, 5/2) -- (0,3);
\draw (-1.7320508,0) -- (-0.8660254,-1/2);
\draw (1.7320508,0) -- (0.8660254,-1/2);
\draw (-0.8660254,-1/2) -- (-0.8660254,-3/2);
\draw (0.8660254,-1/2) -- (0.8660254,-3/2);
\draw (-0.8660254,-3/2) -- (-1.7320508,-2);
\draw (0.8660254,-3/2) -- (1.7320508,-2);
\draw (-1.7320508,0) -- (-2.5980762,-1/2);
\draw (1.7320508,0) -- (2.5980762,-1/2);
\draw (-2.5980762,-3/2) -- (-2.5980762,-1/2);
\draw (2.5980762,-3/2) -- (2.5980762,-1/2);
\draw (-2.5980762,-3/2) -- (-1.7320508,-2);
\draw (2.5980762,-3/2) -- (1.7320508,-2);
\draw (0,-3) -- (-1.7320508,-2);
\draw (0,-3) -- (1.7320508,-2);
\draw (0,-3) -- (0,-6);
\draw (-2.5980762,-1/2) -- (-2.5980762,3/2);
\draw (2.5980762,-1/2) -- (2.5980762,3/2);
\draw (2.5980762,3/2) -- (0.8660254,5/2);
\draw (-2.5980762,3/2) -- (-0.8660254,5/2);
\draw (0, -6) -- (5.1961524, -3);
\draw (0,-6) -- (-5.1961524,-3);
\draw (5.1961524,-3) -- (5.1961524,3);
\draw (-5.1961524,-3) -- (-5.1961524,3);
\draw (0,6) -- (5.1961524,3);
\draw (-5.1961524,3) -- (0,6);
\draw (-5.1961524,3) -- (-2.5980762,3/2);
\draw (5.1961524,3) -- (2.5980762,3/2);
\end{tikzpicture}
\end{subfigure}
\begin{subfigure}{}
\begin{tikzpicture}[scale=0.50]
\draw (0,0) -- (0,1);
\draw (0,0) -- (-0.8660254, -1/2);
\draw (0,0) -- (0.8660254, -1/2);
\draw (0,1) -- (-0.8660254, 3/2);
\draw (0,1) -- (0.8660254, 3/2);
\draw[dashed] (-0.8660254, 3/2) -- (0.8660254, 3/2); 
\draw (-0.8660254, 3/2) node {$\bullet$}; 
\draw (0.8660254, 3/2) node {$\bullet$}; 
\draw (0.8660254, 3/2) -- (0.8660254, 5/2);
\draw (-0.8660254, 3/2) -- (-0.8660254, 5/2);
\draw (0.8660254, 5/2) -- (0,3);
\draw (-0.8660254, 5/2) -- (0,3);
\draw (-1.7320508,0) -- (-0.8660254,-1/2);
\draw (1.7320508,0) -- (0.8660254,-1/2);
\draw[dashed] (-1.7320508,0) to[bend left] (1.7320508,0); 
\draw (-1.7320508,0) node {$\bullet$}; 
\draw (1.7320508,0) node {$\bullet$}; 
\draw (-0.8660254,-1/2) -- (-0.8660254,-3/2);
\draw (0.8660254,-1/2) -- (0.8660254,-3/2);
\draw (-0.8660254,-3/2) -- (-1.7320508,-2);
\draw (0.8660254,-3/2) -- (1.7320508,-2);
\draw[dashed] (-0.8660254,-3/2) -- (0.8660254,-3/2); 
\draw (-0.8660254,-3/2) node {$\bullet$}; 
\draw (0.8660254,-3/2) node {$\bullet$}; 
\draw (-1.7320508,0) -- (-2.5980762,-1/2);
\draw (1.7320508,0) -- (2.5980762,-1/2);
\draw (-2.5980762,-3/2) -- (-2.5980762,-1/2);
\draw (2.5980762,-3/2) -- (2.5980762,-1/2);
\draw (-2.5980762,-3/2) -- (-1.7320508,-2);
\draw (2.5980762,-3/2) -- (1.7320508,-2);
\draw[dashed] (-2.5980762,-3/2) .. controls (-1.7320508,-4) and (1.7320508,-4) .. (2.5980762,-3/2); 
\draw (-2.5980762,-3/2) node {$\bullet$}; 
\draw (2.5980762,-3/2) node {$\bullet$}; 
\draw (0,-3) -- (-1.7320508,-2);
\draw (0,-3) -- (1.7320508,-2);
\draw (0,-3) -- (0,-6);
\draw (-2.5980762,-1/2) -- (-2.5980762,3/2);
\draw (2.5980762,-1/2) -- (2.5980762,3/2);
\draw (2.5980762,3/2) -- (0.8660254,5/2);
\draw (-2.5980762,3/2) -- (-0.8660254,5/2);
\draw (0, -6) -- (5.1961524, -3);
\draw (0,-6) -- (-5.1961524,-3);
\draw (5.1961524,-3) -- (5.1961524,3);
\draw (-5.1961524,-3) -- (-5.1961524,3);
\draw (0,6) -- (5.1961524,3);
\draw (-5.1961524,3) -- (0,6);
\draw (-5.1961524,3) -- (-2.5980762,3/2);
\draw (5.1961524,3) -- (2.5980762,3/2);
\draw[dashed] (-5.1961524,-3) to[bend right] (5.1961524,-3); 
\draw (-5.1961524,-3) node {$\bullet$}; 
\draw (5.1961524,-3) node {$\bullet$}; 
\end{tikzpicture}
\end{subfigure}
\caption{Top: Part of the tree of solutions to $x^2 + y^2 + z^2 = xyz$ in positive integers. Bottom left: The 28 solutions mod 7 connected by Markoff moves $m_1, m_2, m_3$. Bottom right: The same solutions connected by $D_1, D_2, D_3$ where $D_j = \tau_{23} \circ m_j$.
}
\label{fig:mod7}
\end{figure}
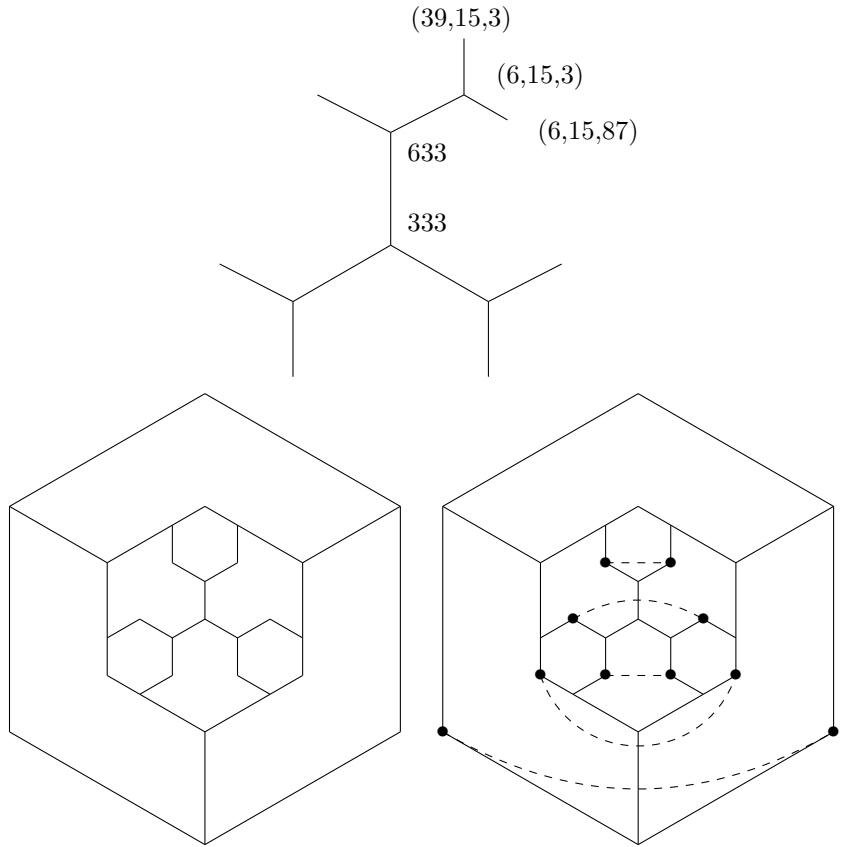

For example, we consider $p=7$ in detail. Using the form $x^2 + y^2 +z^2 = xyz$, we start from the base solution $(3,3,3)$. The Markoff moves lead to $(6,3,3), (3,6,3),$ and $(3,3,6)$. At the second level, one finds the permutations of $(1,6,3) = m_1 (3,6,3)$ because $6 \times 3 - 3 \equiv 1 \bmod 7$. In characteristic 0, we would have $(15, 6, 3)$ instead of $(1,6,3)$, and instead of a loop $m_3(1,6,3) = (1,6,3)$ we would simply have $m_3(15, 6, 3) = (15, 6, 87)$. At the third level come $(3,1,4) = m_3(3,1,6)$ and its six permutations. At the fourth level the relation $(1,1,4)=m_1(3,1,4) = m_2(1,3,4)$ and its permutations lead to three cycles of length 6. The remaining Markoff move leads to, for instance, $(3,4,4) = m_2(3,1,4)$ (and two other permutations). Three cycles of length 8 are formed between $(3,3,3)$ and the permutations of $(3,4,4)$. At the fifth level, one obtains three points of the form $(4,4,6) = m_3(4,4,3)$. At the sixth level, one has all the solutions, with three points of the form $(4,6,6)$ completing a final cycle of length 6 together with the points $(4,4,6)$.
This procedure is shown graphically in Figure~\ref{fig:mod7}.

\begin{figure} [h]
\includegraphics[width=0.25\textwidth]{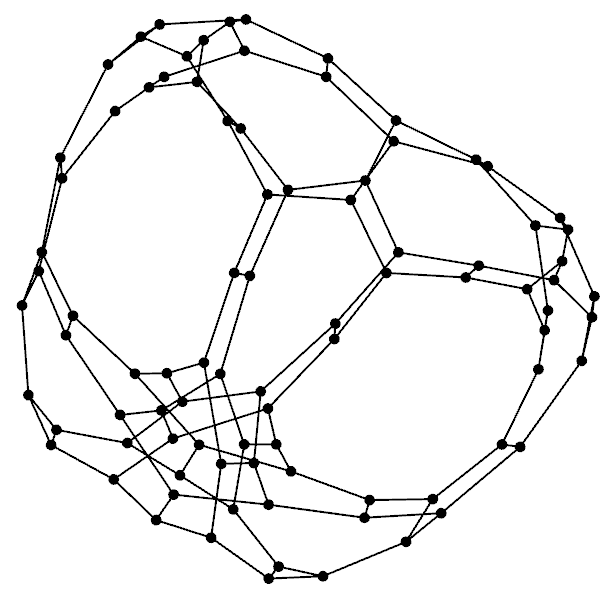}
\caption{Markoff graph for $\F_{11}$}
\label{fig:markoffgraph}
\end{figure}

The Markoff graph mod 7 has 12 loops, at the points $(4,6,6), (1,1,4), (1,6,3)$, and their permutations. If we use the generators $D_1, D_2, D_3$ instead of the Markoff moves, then $(4,6,6)$ and $(4,1,1)$ will still be fixed by $D_1$ but the other points will be joined pairwise. Note that 6 and 1 are the two square roots of $8 \equiv 1 \bmod 7$. Also, the Dehn-neighbours of $(3,3,6)$ are the Markoff-neighbours of $(3,6,3)$. Thus using $D_1, D_2, D_3$ mod 7 has reflected the graph left-to-right around $(3,3,3)$ and created four extra edges.
As another example, the graph constructed from $D_1, D_2, D_3$ for $p=11$ is shown in Figure~\ref{fig:markoffgraph}. It has no loops, since 2 is not a square modulo 11.

The same Vieta operations act on Markoff-like equations $x^2 + y^2 + z^2 - 3xyz = k$ for other values of $k$ besides $0$. The resulting graphs need not be connected, and the structure of the connected components depends on arithmetic relations between $p$ and $k$. 
For example, if $k$ is a square modulo $p$, then there will be a component of size 6 containing $(0,0,\pm\sqrt{k})$ and its permutations. 
In Section \ref{sec:nonzerok}, we give more examples and discuss some patterns in the component sizes (Table \ref{table:components}). 

\begin{table}[h]
\begin{tabular}{c|c|c|c|c}
$k \backslash p$  & 5 & 7 & 11 & 13 \\ \hline
0 & 1 40 & 1 28 & 1 88 & 1 208 \\ \hline
1 & 4 6 16 & 6 16 & 6 160 & 6 112 \\ \hline
2 & 36 & 4 6 16 24 & 144 & 196 \\ \hline
3 & 16 & 64 & 6 160 & 6 216 \\ \hline
4 & 6 & 64 & 6 160 & 6 112 \\ \hline
5 &  & 36 & 6 72 & 144 \\ \hline
6 &  & 64 & 40 60 & 128 16 \\ \hline
7 & & & 144 & 144 \\ \hline
8 & & & 40 60 & 196 \\ \hline
9 & & & 4 6 16 48 48 & 6 216 \\ \hline
10 & & & 16 128 & 6 112 \\ \hline
11 & & & & 196 \\ \hline
12 & & & & 4 6 16 48 96 \\ \hline
\end{tabular}
\caption{Each entry lists the sizes of the orbits of $(x,y,z)$ satisfying $x^2 + y^2 + z^2 - 3xyz = k$ mod $p$ under the group generated by Dehn twists, permutations, and double sign changes. The prime $p$ runs horizontally and the level $k$ runs horizontally. Each column can be extended periodically since $k$ is taken modulo $p$.}
\label{table:components-intro}
\end{table}

Most conspicuously, for each $p > 3$, there is a particular value $k \equiv 4/9 \bmod p$ such that the graph has a large number of orbits compared to other level sets. In the normalization $x^2+y^2+z^2=xyz + k$ instead of $3xyz$, this special value is simply $k=4$. This degenerate level set is called the \emph{Cayley cubic}:
\begin{equation} \label{eqn:cayley-cubic}
x^2 + y^2 + z^2 = xyz + 4.
\end{equation}
In Section~\ref{sec:cayley}, we discuss this case in detail and prove the following theorem.
\begin{theorem} \label{thm:main}
For each prime $p \geq 5$, the orbits of the action of Markoff moves and permutations on solutions of (\ref{eqn:cayley-cubic}) modulo $p$ are in bijection with those divisors of $p^2-1$ that are multiples of $p+1$ or $p-1$.
The number of orbits is
\begin{equation} \label{eqn:num-orb}
e_2  \prod_{q \in Q^-} (e_q+1) + 2 \prod_{q \in Q^+} (e_q+1)  - 2.
\end{equation}
where $Q^+$ and $Q^-$ are disjoint sets of odd primes dividing $p^2-1$ such that the prime factorizations of $p+1$ and $p-1$ are
\[
p + \varepsilon = 2 \prod_{q \in Q^+} q^{e_q}, \quad p - \varepsilon = 2^{e_2-1} \prod_{q \in Q^-} q^{e_q}.
\]
We write $\varepsilon = (-1)^{\frac{p-1}{2}}$.
Given such a divisor $\prod_q q^{f_q} \neq p^2-1, \frac{p^2-1}{2}$, the size of the corresponding orbit is
\begin{equation} \label{eqn:size-orb}
\frac{1}{2} \prod_{q: f_q < e_q} (q^2-1) q^{2(e_q - f_q - 1)}
\end{equation}
whereas $p^2-1$ and $\frac{p^2-1}{2}$ correspond to orbits of twice this size.
\end{theorem}

We illustrate the notation of Theorem~\ref{thm:main} for $p=5, 7, 11$ in Table~\ref{table:orbits}, and give some larger examples in Section~\ref{sec:cayley}. 
Some of the orbits described in Theorem~\ref{thm:main} merge when we include the further symmetries
$(x, y, z) \mapsto (\sigma_1 x, \sigma_2 y, \sigma_3 z)$ with signs obeying  $\sigma_1 \sigma_2 \sigma_3 = 1$. 
In terms of divisors, we show in Section~\ref{sec:cayley} that the effect is to join the orbits of $t$ and $2t$ when the power of 2 dividing $t$ is just one less than the power dividing $p^2-1$.
Thus we have the following corollary, also illustrated in Table~\ref{table:orbits}.
\begin{corollary} \label{cor:signs}
The number of orbits in the Cayley cubic, including the effect of sign changes, is
\begin{equation} \label{eqn:num-orb-sign}
(e_2 -1 ) \prod_{q \in Q^-} (e_q+1) +  \prod_{q \in Q^+} (e_q+1) - 1.
\end{equation}
\end{corollary}

Another consequence of Theorem~\ref{thm:main} is that the number of orbits is small compared to $p$.
\begin{corollary} \label{cor:upper}
The number of orbits for the Cayley cubic modulo $p$ is at most the number of divisors of $p^2-1$. Hence for any $\delta > 0$ there is a number $A_{\delta} > 0$ such that there are at most $A_{\delta} p^{\delta}$ orbits.
\end{corollary}

On the other hand, (\ref{eqn:num-orb}) shows that the number of orbits is at least $e_2$.
Let $p$ be congruent to $\pm 1 \bmod 2^k$ with $k \geq 2$. 
Dirichlet's theorem guarantees that there are infinitely many such primes for each $k$. The first one is at most $2^{kL}$ where $L$ is Linnik's constant.
For such primes,
\[
e_2 = k+1 > \frac{1}{L} \log_2{p}
\]
so Theorem~\ref{thm:main} implies that along an infinite subsequence beginning $p = 5, 7, 17, 23, \ldots$ of primes congruent to $\pm 1 \bmod 2^k$ for $k = 1, 2, 3, \ldots$, the number of orbits grows at least logarithmically.
See Linnik's original articles \cite{L1}, \cite{L2} proving that there is a finite such $L$, and \cite{X} for a recent numerical value $L \leq 5$ obtained by Xylouris. Assuming the Riemann Hypothesis, one could take $L$ arbitrarily close to 2.
If $p = 2^l - 1$ is a Mersenne prime, then
\[
e_2 = l+1 = \log_2(p+1)+1
\]
and it follows that the Cayley cubic modulo $p$ has more than $\log_2{p}$ orbits in these cases.
Over all primes up to a given magnitude, the average number of orbits is of order $\log{p}$, along the same lines as Titchmarsh's divisor problem \cite{T}.

It may be surprising that the factors of $p^2-1$ prove decisive for the orbit structure modulo $p$. A preliminary change of variable naturally leads one to work in the extension $\F_{p^2}$ rather than $\F_p$, and $p^2-1$ is the order of the multiplicative group $\F_{p^2}^{\times}$. The special feature of (\ref{eqn:cayley}) is that the action linearizes, and our approach in Section~\ref{sec:cayley} is to determine the orbits explicitly by ``linear algebra mod $p^2-1$".

\begin{table}[h]
\begin{tabular}{clclclclclclclclc}
$p$ & $\varepsilon$ & $p^2-1$ & $Q^+$   & $Q^-$       & $e_2$ & Orbit sizes & $\pm$-orbit sizes \\ \hline
5   & 1          & 24          & $\{3\}$ & $\emptyset$ & 3     & 1 3 4 6 12                  & 4 6 16                      \\
7   & $-1$         & 48          & $\{3\}$ & $\emptyset$ & 4     & 1 3 4 6 12 24               & 4 6 16 24                   \\
11  & $-1$         & 120         & $\{5\}$ & $\{3\}$     & 3     & 1 3 4 6 12 12 36 48         & 4 6 16 48 48               
\end{tabular}
\caption{Examples of the notation in Theorem~\ref{thm:main}, showing the sizes of the orbits under permutations and Markoff moves as well as the sizes of the orbits after including double sign changes.}
\label{table:orbits}
\end{table}

We recommend the book by Aigner \cite{A} and the article by Bombieri \cite{Bo} for more information about the Markoff equation over $\Z$ and its many interconnections with other subjects.
It is known roughly how many solutions there are subject to a given bound on the coordinates. Zagier proved that the number of solutions with $x,y$, and $z$ less than $T$ is asymptotic to $C\log(T)^2$, with $C = 0.1807...$ given by an explicit infinite series \cite{Z}. Mirzakhani gave a proof making use of the Markoff surface's beautiful relation to trace identities and geometry \cite{Mirz}. We review  Fricke's trace identity in Section~\ref{sec:cayley}. Because of the connection between the Markoff surface and trace identities, the connectedness (or not) of different level sets is related to group-theoretic conjectures of McCullough and Wanderley \cite{MW}.
We begin more modestly with a formula (Proposition \ref{prop:numel_levelset}) for the number of solutions mod $p$ to the Markoff equation (\ref{eqn:markoff}) and equations of a similar form. 

\section{Number of solutions} \label{sec:numel}

For a finite field $\F_p$, simply counting the triples $(x,y,z)$ with no regard for whether they solve (\ref{eqn:markoff}) or not shows that the Markoff equation has at most $p^3$ solutions, and one expects approximately $p^2$ solutions since a single equation is imposed. Using the Legendre symbol to detect how many roots a quadratic equation has, we give an exact formula for the number of solutions to (\ref{eqn:markoff}) modulo any prime $p \geq 5$. This is the number of vertices in the graph, which is useful to know in advance when we come to enumerate the connected components.

\begin{proposition} \label{prop:numel_levelset}
The number of solutions to the equation $x^2+y^2+z^2-axyz=k$ with $a\neq 0$ in a finite field $\F_p$ is
\begin{equation}
N(p,k,a) = p^2 + \left(3 + \legendre{k}{p}\right)\legendre{a^2k - 4}{p} p + 1.
\end{equation}
In particular, for the Markoff equation and its level sets $x^2 + y^2 + z^2 - 3xyz - k = 0$, the number of solutions mod $p$ is
\begin{equation}
N(p,k)=
\begin{cases}
p^2 + 3\legendre{-1}{p}p + 1 & \text{if } k = 0 \\
p^2 + 1 & \text{if } 9k = 4 \\
p^2 + 4\legendre{9k-4}{p}p + 1 & \text{otherwise} 
\end{cases}
\end{equation}
\end{proposition}
The special case of the original Markoff surface -- $a=3$ and $k=0$ in the notation above -- was previously calculated by Baragar in \cite{Bar} pages 117-122, which in this case also gives the number of solutions over fields with $p^k$ elements rather than just the prime fields. See also the note of Carlitz \cite{C} for the general case.

Under a scaling $(x,y,z) = b \cdot (\xi, \eta, \zeta)$, the equation $x^2 + y^2 + z^2 - axyz = k$ transforms into $\xi^2 + \eta^2 + \zeta^2 - ab\xi \eta \zeta = k/b^2$. Thus we can fix a convenient value such as $a = 1$ or $a = 3$ with no loss of generality.
The degenerate case where $a^2 k = 4$ comes from the \emph{Cayley cubic}, which we will see leads to a distinctive graph structure.

\begin{proof}

We sum over all $x$ and $y$ in $\F_p$, tallying how many $z$ satisfy equation (\ref{eqn:markoff}) for each pair $(x,y)$. There are 2, 1, or 0 solutions modulo $p$ according to whether the discriminant of the resulting quadratic is a quadratic residue, 0, or a nonresidue. We can therefore express $N(p,k,a)$ using the Legendre symbol:
\begin{equation}
N(p,k,a) = \sum_{x \in \F_p} \sum_{y \in \F_p} \left(1 + \legendre{(axy)^2 - 4(x^2 + y^2 - k)}{p} \right).
\end{equation}
Summing the term 1 over all pairs $(x,y)$ yields $p^2$, which is the main contribution to $N(p,k,a)$ claimed in Proposition~\ref{prop:numel_levelset}. 
Now we evaluate the lower-order contribution from the character sum:
\begin{equation*}
S =  \sum_{x \in \F_p} \sum_{y \in \F_p}  \legendre{(axy)^2 - 4(x^2 + y^2 - k)}{p} .
\end{equation*}
Fixing $x$, we sum the Legendre symbol of $y^2 ((ax)^2 - 4) + 4(k-x^2)$ over $y$. For two values of $x$, namely $\pm 2/a$, we have $(ax)^2 - 4 = 0$ and then all $p$ of the inner summands are
\begin{equation*}
\legendre{4(k-x^2)}{p} = \legendre{k-x^2}{p} = \legendre{k - 4/a^2}{p} = \legendre{a^2k - 4}{p}.
\end{equation*}
Thus
\begin{equation*}
S = 2 \legendre{a^2 k - 4}{p}p + \sum_{(ax)^2 \neq 4} \legendre{(ax)^2-4}{p} \sum_y \legendre{y^2 + 4(k-x^2)/((ax)^2 - 4)}{p}
\end{equation*}
The inner sum over $y$ is a sum of Legendre symbols of \emph{shifted squares} $y^2 - s$. If $s=0$, every term is 1, except for the term 0 when $y=0$. Thus the sum is $p-1$ in case $k = x^2$. Otherwise, we use the following Lemma:
\begin{lemma}
\label{lem:pedro}
For any nonzero $s \in \F_p$,
$$
\sum_{z \in \F_p} \legendre{z^2-s}{p} = -1.
$$
\end{lemma}
Let us finish the calculation of $S$ and then prove Lemma~\ref{lem:pedro}. By the lemma, the sum over $y$ is $-1$, unless $k = x^2$ in which case it is $p-1$. The contribution from $x^2 = k$ is
\begin{equation*}
(p-1)\legendre{a^2 k - 4}{p} \left(1 + \legendre{k}{p} \right)
\end{equation*}
which, in particular, is 0 if there are no such $x$. The remaining contribution to $S$ is the sum
\begin{equation*}
\sum_{k \neq x^2 \neq 4/a^2} \legendre{(ax)^2 - 4}{p} (-1).
\end{equation*}
The constraint $x^2 \neq 4/a^2$ may be ignored since the summand is 0 in that case. We have another sum over shifted squares $u^2 - s$ with $u = ax$ and $s = 4$, except with solutions to $u^2 = a^2 k $ omitted if there are any to begin with. Therefore, by Lemma~\ref{lem:pedro} again,
\begin{equation*}
\sum_{k \neq x^2 \neq 4/a^2} \legendre{a^2 x^2 - 4}{p} = -1 - \legendre{a^2k - 4}{p}\left(1 + \legendre{k}{p}\right).
\end{equation*}
Including all the cases, we have
\begin{align*}
S &= 2\legendre{a^2k-4}{p} p + \left(1 + \legendre{k}{p}\right)\legendre{a^2k - 4}{p} (p-1) - \left(-1 - \legendre{a^2k - 4}{p}\left(1 + \legendre{k}{p}\right) \right) \\
&= \legendre{a^2k - 4}{p} \left( 3 + \legendre{k}{p} \right) p + 1,
\end{align*}
which implies that $N(p,k,a) = p^2 + S$ is as claimed.
\end{proof}

To prove Lemma~\ref{lem:pedro}, first consider the case that $s$ is a quadratic residue, say $s = t^2$ with $t \neq 0$. Then
\begin{align*}
\sum_{z \in \F_p} \legendre{z^2 - s}{p} &= \sum_z \legendre{z-t}{p} \legendre{z + t}{p} \\
&= 0 + \sum_{u \neq 0} \legendre{u}{p} \legendre{u + 2t}{p} 
\end{align*}
having changed variables to $u = z-t$, which runs over $\F_p$ just as $z$ does but gives no contribution when $u=0$.
The inverse $u^{-1}$ is a quadratic residue precisely when $u$ is, so
\[
\legendre{u}{p} \legendre{u+2t}{p} = \legendre{u^{-1}}{p} \legendre{u+2t}{p} = \legendre{1+2t/u}{p}.
\]
As $u$ ranges over all non-zero values, $1 + 2t/u$ ranges over all values except 1. The sum of all the Legendre symbols is 0, so when we omit $1 = \legendre{1}{p}$, the sum is $-1$. This proves the lemma in case $s$ is a quadratic residue. Next observe that the sum in Lemma~\ref{lem:pedro} only depends on whether $s$ is a quadratic residue. Thus it is $p-1$ when $s = 0$, $-1$ when $s$ is a quadratic residue, and some value $n$ when $s$ is a quadratic non-residue. On the other hand
\begin{equation*}
\sum_s \sum_z \legendre{z^2 - s}{p} = \sum_z \sum_s \legendre{z^2 - s}{p} = 0
\end{equation*}
which implies that $p-1 - (p-1)/2 + n(p-1)/2 = 0$, so $n$ must also be $-1$.

\section{Connectedness and Spectral Gap for $k=0$} \label{sec:zerok}

The Cheeger constant $h(G)$ measures whether there is a bottleneck in the graph $G$. It is defined as
\begin{equation} \label{eqn:cheegerdefn}
h(G) = \min \frac{|\partial A|}{|A|}
\end{equation}
where $A$ is any nonempty subset of $G$ that has at most half the vertices of $G$ and $\partial A$ is its edge boundary. If $h(G) = 0$, then $G$ is disconnected since there is a subset $A$ with $|\partial A| = 0$, that is, no edges from $A$ to $G \setminus A$. A large value of $h(G)$ means that there are many ways to escape from any given $A$. This is referred to as \emph{expansion}. For a graph $G$ with $V$ vertices, there are about $2^{V-1}$ subsets of $G$ containing at most half the vertices of $G$, so there are an exponential number of candidates for the minimum in (\ref{eqn:cheegerdefn}). Therefore, instead of computing the exact value of $h(G)$, it is practical to estimate it. The \emph{Cheeger inequality} states that
\begin{equation} \label{eqn:cheegerineq}
\frac{1}{2}(d-\lambda_2) \leq h(G) \leq \sqrt{2d(d-\lambda_2)}
\end{equation}
where $G$ is a $d$-regular graph and $\lambda_2$ is the second highest eigenvalue of the adjacency matrix of $G$ after $d$. It is named after an analogous theorem of Cheeger on manifolds, the theorem for graphs being due to Dodziuk and Alon-Milman in \cite{D}, \cite{AM}. The difference $d-\lambda_2$ is known as the \emph{spectral gap} of $G$. Note that $d$ is always an eigenvalue of a $d$-regular graph because the ``all 1's vector" is an eigenvector when all the rows sum to $d$.

In particular, (\ref{eqn:cheegerineq}) shows that $h(G)$ and $d - \lambda_2$ converge to 0 or not together. Therefore, the spectral gap is an equally good measure of how well connected a graph is. A lower bound on $d - \lambda_2$ indicates the extent to which the graph is well connected. The advantage of the spectral notion is that $\lambda_2$ is easier to compute than $h(G)$. The Markoff equation in $\F_p$ has roughly $p^2$ solutions, so $G$ is represented by a $p^2$-by-$p^2$ matrix, and $\lambda_2$ is its second highest eigenvalue.
Although computing the eigenvalues of such an enormous matrix is costly, it is much better than checking the roughly $2^V \approx 2^{p^2}$ subsets $A$ needed to compute $h(G)$ by brute force.

\begin{figure}[h] 
\centering
\includegraphics[width=\textwidth]{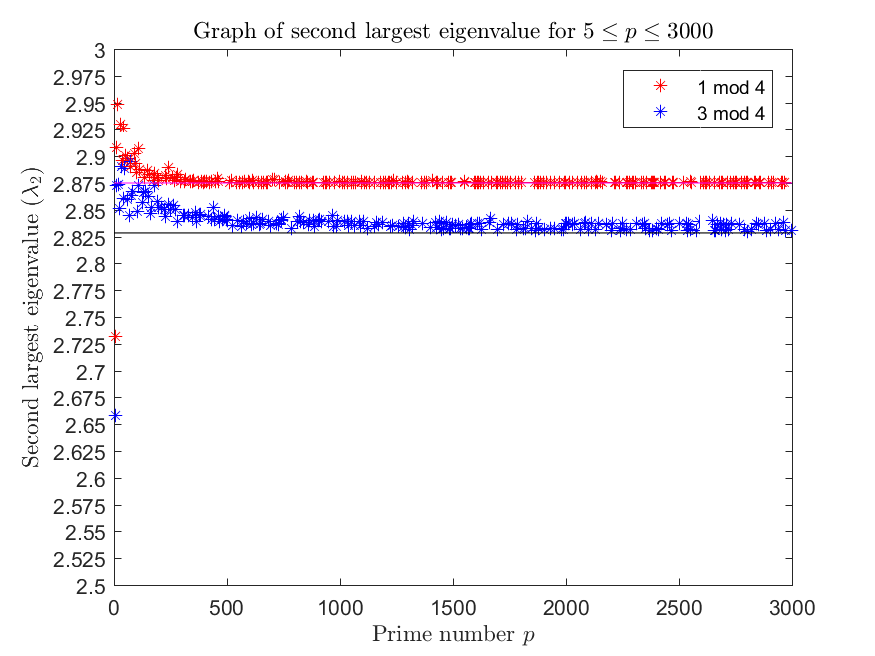}
\caption{Graph of $\lambda_2$ for $5 \leq p \leq 3000$. The black line marks the Ramanujan case $\lambda_2 = 2\sqrt{2} = 2.828...$}
\label{fig:eigenvalue}
\end{figure}

Figure~\ref{fig:eigenvalue} shows a striking pattern in the values of $\lambda_2$ for primes less than 3000. The black horizontal line marks $2\sqrt{2}=2.828\ldots$, and the magenta line marks $2.875$. For Markoff graphs modulo a prime $p$ congruent to 3 modulo 4, the data suggests that $\lambda_2$ approaches $2\sqrt{2}$. For prime numbers $p$ congruent to 1 modulo 4, the data suggests that $\lambda_2$ approaches a higher value. Thus for primes congruent to 1 modulo 4, the Markoff graphs seem to exhibit weaker expansion compared to primes congruent to 3 modulo 4.

The apparent limit $2\sqrt{2}$ is a familiar number for 3-regular graphs. A $d$-regular graph is a \emph{Ramanujan graph} if $\lambda_2 \leq 2\sqrt{d-1}$. These graphs are the optimal expanders. See \cite{LPS} for the first construction of such graphs and more information.

\begin{figure}[h] 
\centering
\hfill
\subfigure[$p=83$]{\includegraphics[width=0.48\textwidth]{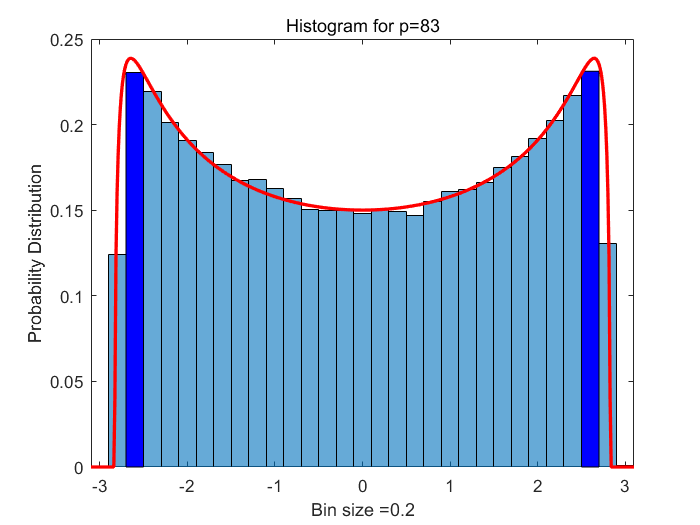}}
\hfill
\subfigure[$p=89$]{\includegraphics[width=0.48\textwidth]{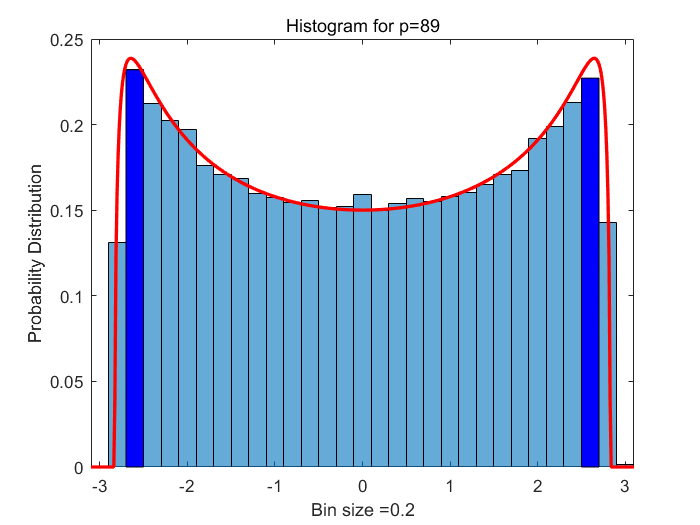}}
\hfill
\caption{Histogram of eigenvalues for $p=83$ and $89$}
\label{fig:histogram}
\end{figure}

\begin{figure}[h] 
\centering
\includegraphics[width=\textwidth]{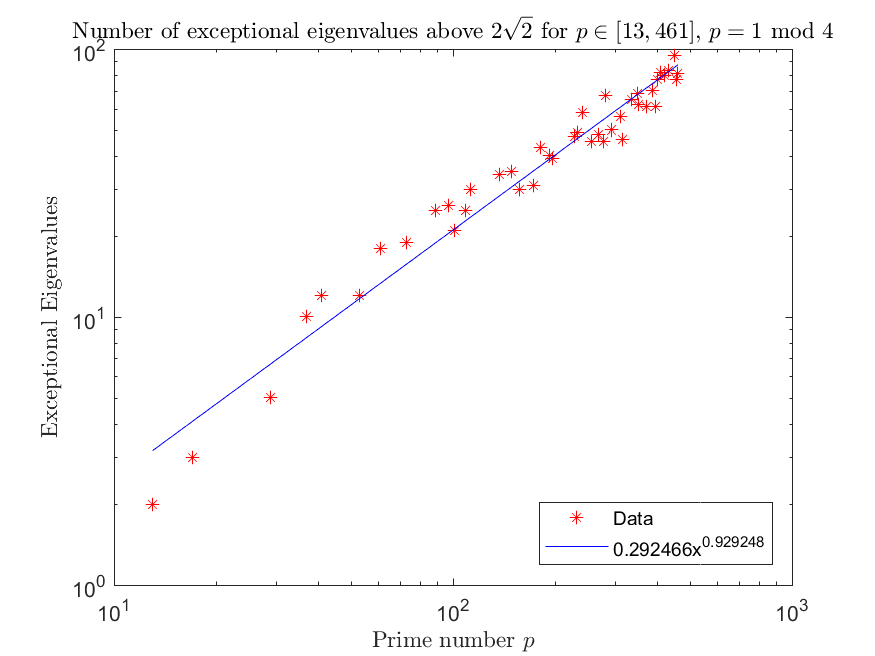}
\caption{Log-Log Plot of Eigenvalues above $2 \sqrt{2}$ for $p$ congruent to 1 modulo 4}
\label{fig:outstanding_eigen}
\end{figure}

Beyond $\lambda_2$, we computed all of the eigenvalues of these matrices for a smaller range of primes. For comparison, the Kesten-McKay Law specifies the eigenvalue distribution of a random $d$-regular graph \cite{Mc}, \cite{K}. It is given by the probability density function
\begin{equation}
\rho_d(\lambda) = \frac{d}{2\pi} \frac{ \sqrt{4(d-1) - \lambda^2} }{d^2 - \lambda^2} \mathbbm{1}_{[-2\sqrt{d-1},2\sqrt{d-1}]}(\lambda)
\end{equation}
In particular, it is supported on the interval $[-2\sqrt{d-1}, 2\sqrt{d-1}]$. For $3$-regular graphs, the distribution is bimodal (the maxima are at $\pm \sqrt{7}$) and supported on the interval $[-2\sqrt{2}, 2\sqrt{2}]$. For both $p$ congruent to 1 modulo 4 and $3$ modulo $4$ alike, the histogram of eigenvalues follows the Kesten-McKay Law closely (Figure~\ref{fig:histogram}). This suggests that although $\lambda_2$ converges to a higher value for $p$ congruent to $1$ modulo $4$, this is only because of a vanishing proportion of exceptional eigenvalues above $2\sqrt{2}$. Indeed, Figure~\ref{fig:outstanding_eigen} seems to indicate that the number of eigenvalues above $2\sqrt{2}$ grows only like $p$ out of the total of roughly $p^2$ eigenvalues.
The Kesten-McKay law for Markoff graphs has recently been proved \cite{DM}, although the resulting bound for the number of exceptional eigenvalues is only $p^2/\log{p}$ instead of $p$.

\section{The Structure of graphs from non-zero $k$} \label{sec:nonzerok}
A variation of the Markoff surface can be created by adding a constant $k$:
\begin{equation} \label{eqn:levelk}
x^2+y^2+z^2-3xyz=k.
\end{equation}
This equation is invariant under the same Vieta moves and permutations of $(x,y,z)$, so the Dehn twists $D_1, D_2, D_3$ act on it exactly as in the case $k=0$. However, for nonzero $k$, the connectedness of the resulting graph is no longer guaranteed. For example, when $p=7$ and $k=2$, there are 10 connected components. Also, the triple $(1,1,1)$ is no longer a guaranteed solution. We instead use a brute-force method to find a solution to serve as the root from which to explore using the generators $D_j$. Eventually, the component of that solution is fully unveiled. If its size agrees with the total number of solutions predicted in Proposition~\ref{prop:numel_levelset}, then we have finished constructing the graph. If not, we must use brute force again to find a solution outside this component and continue exploring from there.

\begin{figure}[h]
\centering
\hfill
\subfigure[]{
\includegraphics[scale=0.25]{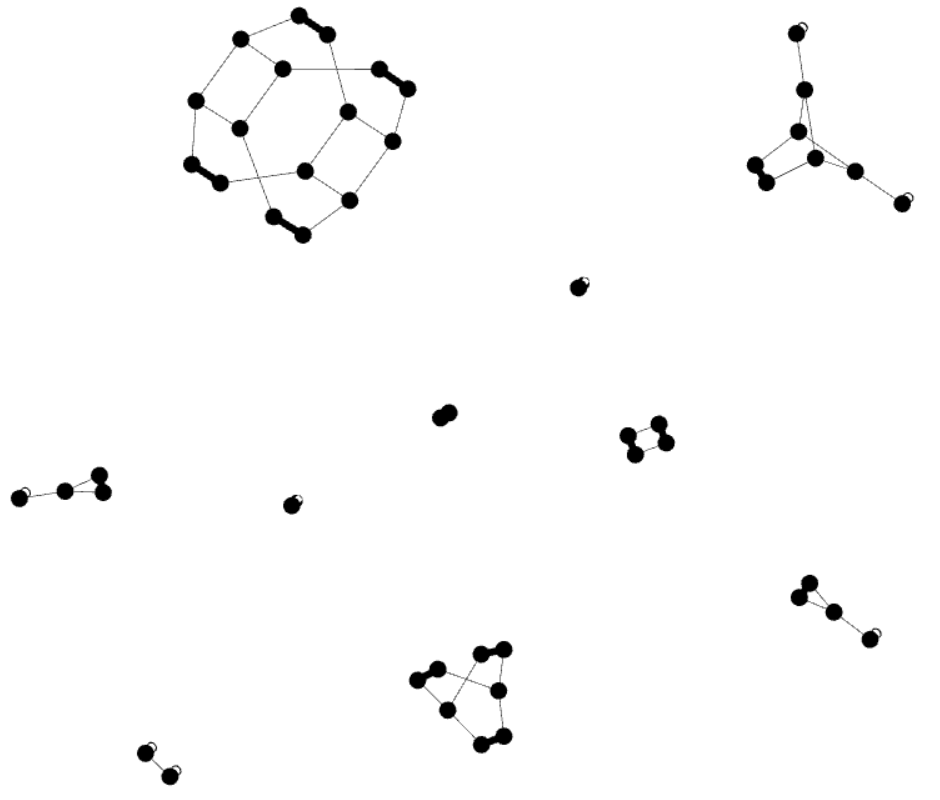}
}
\hfill
\hfill
\hfill
\hfill
\hfill
\hfill
\hfill
\hfill
\subfigure[]{
\includegraphics[scale=0.25]{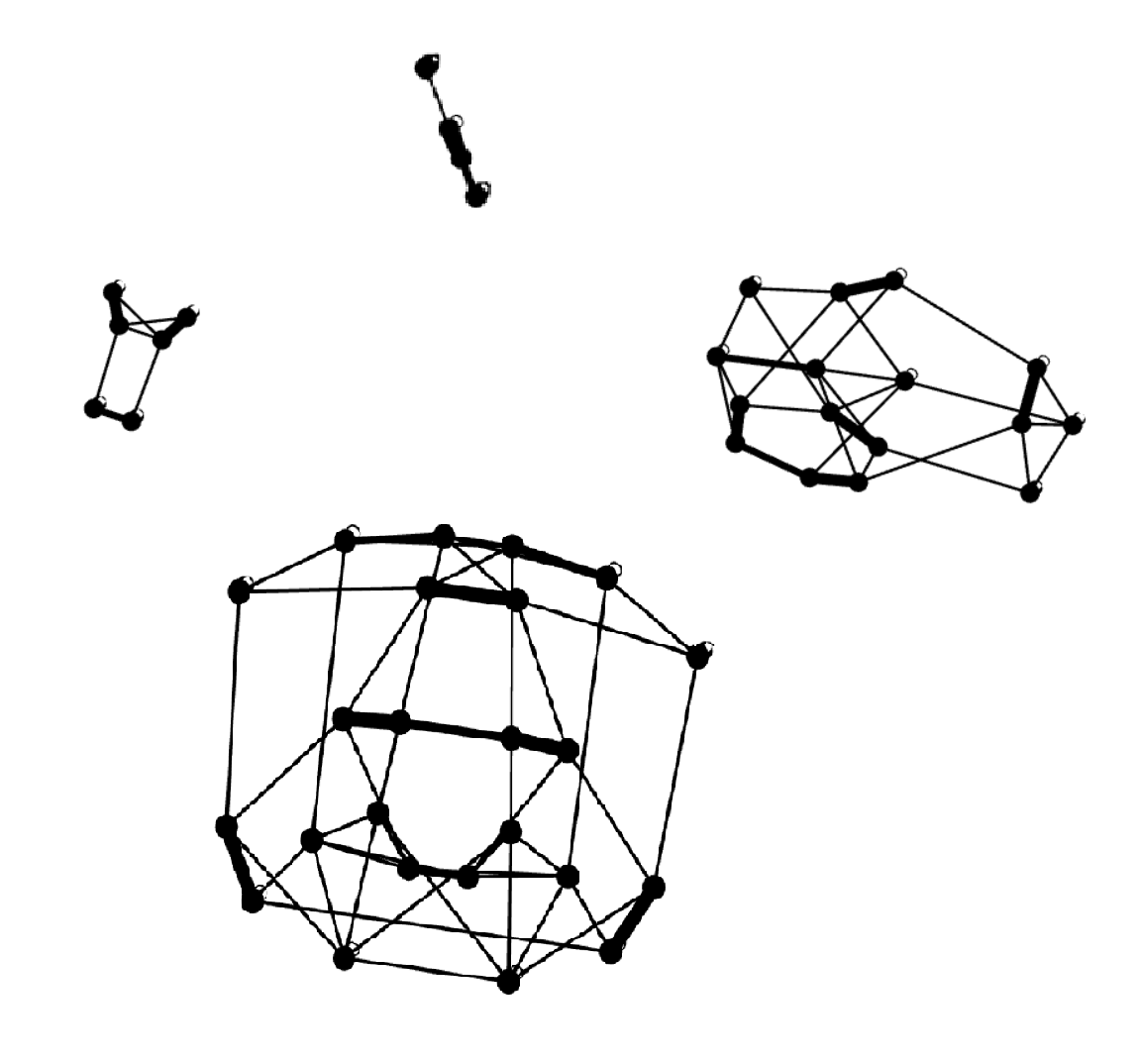}
}
\caption{The solutions to $x^2 + y^2 + z^2 - 3xyz = 2$ mod $7$, with edges corresponding to $D_1,D_2,D_3$ in (a), and then in (b) with additional edges for the transpositions $T_{12}$ and $T_{23}$ and the sign change $N_{12}$. The resulting components have sizes 1, 3, 4, 6, 12, 24 in (a) or 4, 6, 16, and 24 in (b).}
\label{fig:merged_disconnected_graph}
\end{figure}

In an effort to make the Markoff graph connected, we enlarged the generating set by including permutations and double-sign-changes that would connect related vertices. The new operations are
\begin{align*} \label{perm-sign}
T_{12}: & (x,y,z) \mapsto (y,x,z) \\
T_{23}: & (x,y,z) \mapsto (x,z,y) \\
N_{12}: & (x,y,z) \mapsto (-x, -y, z)
\end{align*}
These operations do indeed join some components together, but for many pairs $(p,k)$, even the extended graph is still disconnected. Figure~\ref{fig:merged_disconnected_graph} shows both graphs in the case $p=7$ and $k=2$.

\begin{table}
\begin{tabular}{c|c|c|c|c|c}
$k \backslash p$  & 5 & 7 & 11 & 13 & 17 \\ \hline
0 & 1 40 & 1 28 & 1 88 & 1 208 & 1 340 \\ \hline
1 & 4 6 16 & 6 16 & 6 160 & 6 112 & 6 216 \\ \hline
2 & 36 & 4 6 16 24 & 144 & 196 & 6 216 \\ \hline
3 & 16 & 64 & 6 160 & 6 216 & 256\\ \hline
4 & 6 & 64 & 6 160 & 6 112 & 6 16 336 \\ \hline
5 &  & 36 & 6 72 & 144 & 256 \\ \hline
6 &  & 64 & 40 60 & 128 16 & 36 288\\ \hline
7 & & & 144 & 144 & 324 \\ \hline
8 & & & 40 60 & 196 & 4 6 16 24 96 144 \\ \hline
9 & & & 4 6 16 48 48 & 6 216 & 6 352 \\ \hline
10 & & & 16 128 & 6 112 & 324 \\ \hline
11 & & & & 196 & 256 \\ \hline
12 & & & & 4 6 16 48 96 & 324 \\ \hline
13 & & & & & 6 216 \\ \hline
14 & & & & & 256 \\ \hline
15 & & & & & 6 216 \\ \hline
16 & & & & & 6 352  \\ 
\end{tabular}
\caption{Each entry lists the sizes of the orbits of $(x,y,z)$ satisfying $x^2 + y^2 + z^2 - 3xyz = k$ mod $p$ under the group generated by Dehn twists, permutations, and double sign changes. The prime $p$ runs horizontally and the level $k$ runs horizontally. Each column can be extended periodically since $k$ is taken modulo $p$.}
\label{table:components}
\end{table}

Table~\ref{table:components} shows many patterns. There are a few component sizes that appear many times in the table. For example, when $k$ is a square modulo $p$, there is always a size-6 component. This is a result of permutations of $(0,0,\pm \sqrt{k})$. In particular, the graph is always disconnected when $k=1$. The size-1 component that appears when $k=0$ is the $(0,0,0)$ component which was disregarded during the discussion in Section~\ref{sec:zerok}.

For small primes and $k=1$, we have just a single component besides the size 6 component containing $(0,0,1)$. However, when $p = 41$, there are three components of respective sizes 6, 40, and 1800. This extra component of size 40 stems from the fact that 5 is a square mod 41, so that a finite orbit constructed in characteristic 0 from the golden ratio $(1+\sqrt{5})/2$ appears. We refer to Dubrovin for these characteristic 0 orbits in the context of the braid group, p.244 of \cite{Dub}. Modulo other primes for which 5 is a quadratic residue, there is also a component of size 40 but for different levels rather than $k=1$.

There is one $k$ for each $p$ that generates a Markoff graph with an especially large number of components. In Table~\ref{table:components}, these pairs $(p,k)$ are $(5,1), (7,2), (11,9), (13,12) \dots$. These occur whenever 
\begin{equation}
9k-4 \equiv 0 \ \text{mod} \ p.
\end{equation}
namely $k = (2/3)^2$, with the division by 3 understood modulo $p \geq 5$. For this special value of $k$, the Markoff equation becomes a form of the Cayley cubic surface. In this surface, the operations that generate the graph linearize, which leads to more components than in other cases.

\section{The Cayley Cubic} \label{sec:cayley}
The Cayley cubic is a special (degenerate) cubic surface given by
\begin{equation} \label{eqn:cayley}
x^2 + y^2 + z^2 - xyz = 4.
\end{equation}
This is a special case of the Markoff level set, so the same Dehn twists, permutations, and double sign changes act on its solutions modulo any prime $p$. The sizes of the resulting components are listed in Table \ref{table:cayley-components}.

\begin{table}
\begin{tabular}{c|c|c}
$p$ & Factors of $p^2-1$ & {\rm Component sizes} \\ \hline
5 & $2^3 \cdot 3$ & 4 6 16 \\ \hline
7 & $2^4 \cdot 3$ & 4 6 16 24 \\ \hline
11 & $2^3 \cdot 3 \cdot 5$ & 4 6 16 48 48 \\ \hline
13 & $2^3 \cdot 3 \cdot 7$ & 4 6 16 48 96 \\ \hline
17 & $2^5 \cdot 3^2$ & 4 6 16 24 96 144 \\ \hline
19 & $2^3 \cdot 3^2 \cdot 5$ & 4 6 16 48 144 144 \\ \hline
23 & $2^4 \cdot 3 \cdot 11$ & 4 6 16 24 48 192 240 \\ \hline
29 & $2^3 \cdot 3 \cdot 5 \cdot 7$ & 4 6 16 48 96 288 384 \\ \hline
31 & $2^6 \cdot 3 \cdot 5$ & 4 6 16 24 48 96 384 384 \\ \hline
37 & $2^3 \cdot 3^2 \cdot 19$ & 4 6 16 48 144 432 720 \\ \hline
41 & $2^4 \cdot 3 \cdot 5 \cdot 7 $ & 4 6 16 24 48 96 144 576 768 \\ \hline
43 & $2^3 \cdot 3 \cdot 7 \cdot 11$ & 4 6 16 96 240 720 768 \\ \hline
47 & $2^5 \cdot 3 \cdot 23$ & 4 6 16 24 48 96 192 768 1056 \\ \hline
53 & $2^3 \cdot 3^3 \cdot 13$ & 4 6 16 144 336 1008 1296 \\ \hline
59 & $2^3 \cdot 3^3 \cdot 13 $ & 4 6 16 48 48 144 384 1152 1680 \\ \hline
61 & $2^3 \cdot 3 \cdot 5 \cdot 31$ & 4 6 16 48 48 144 384 1152 1920 \\ \hline
67 & $2^3 \cdot 3 \cdot 11 \cdot 17$ & 4 6 16 240 576 1728 1920 \\ \hline
71 & $2^4 \cdot 3^2 \cdot 5 \cdot 7 $ & 4 6 16 24 48 48 96 144 192 432 1728 2304 \\ \hline
73 & $2^4 \cdot 3^2 \cdot 37$ & 4 6 16 24 48 144 192 432 1728 2736 \\ \hline
79 & $2^5 \cdot 3 \cdot 5 \cdot 13$ & 4 6 16 24 48 96 144 336 576 2304 2688 \\ \hline
83 & $2^3 \cdot 3 \cdot 7 \cdot 41$ & 4 6 16 48 96 288 768 2304 3360 \\ \hline
89 & $2^4 \cdot 3^2 \cdot 5 \cdot 11$ & 4 6 16 24 48 144 240 384 720 2880 3456 \\ \hline
97 & $2^6 \cdot 3 \cdot 7^2$ & 4 6 16 24 48 96 96 192 384 768 3072 4704 \\ \hline
101 & $2^3 \cdot 3 \cdot 5^2 \cdot 17$ & 4 6 16 48 144 576 1200 3600 4608 \\ \hline
103 & $2^4 \cdot 3 \cdot 13 \cdot 17$ & 4 6 16 24 336 576 1008 4032 4608 \\ \hline
107 & $2^3 \cdot 3^3 \cdot 53$ & 4 6 16 48 144 432 1296 3888 5616 \\ \hline
109 & $2^3 \cdot 3^3 \cdot 5 \cdot 11$ & 4 6 16 48 48 144 240 432 1296 3888 5760 \\ \hline
113 & $2^5 \cdot 3 \cdot 7 \cdot 19$ & 4 6 16 24 96 96 288 720 1152 4608 5760 \\ \hline
127 & $2^8 \cdot 3^2 \cdot 7$ & 4 6 16 24 96 96 144 384 768 1536 6144 6912 \\ \hline
131 & $2^3 \cdot 3 \cdot 5 \cdot 11 \cdot 13$ & 4 6 16 48 48 240 336 720 1920 5760 8064 \\ \hline
137 & $2^4 \cdot 3 \cdot 17 \cdot 23$ & 4 6 16 24 576 1056 1728 6912 8448 \\ \hline
139 & $2^3 \cdot 3 \cdot 5 \cdot 7 \cdot 23$ & 4 6 16 48 96 144 288 1056 2304 6912 8448 \\ \hline
149 & $2^3 \cdot 3 \cdot 5^2 \cdot 37$ & 4 6 16 48 384 1200 2736 8208 9600 \\ \hline
151 & $2^4 \cdot 3 \cdot 5^2 \cdot 19$ & 4 6 16 24 48 384 720 1200 2160 8640 9600 \\ \hline
157 & $2^3 \cdot 3 \cdot 13 \cdot 79$ & 4 6 16 48 336 1008 2688 8064 12480 \\ \hline
163 & $2^3 \cdot 3^4 \cdot 41$ & 4 6 16 144 1296 3360 10080 11664 \\ \hline
167 & $2^4 \cdot 3 \cdot 7 \cdot 83$ & 4 6 16 24 48 96 192 288 768 1152 2304 9216 13776 \\ \hline
173 & $2^3 \cdot 3 \cdot 29 \cdot 43$ & 4 6 16 1680 3696 11088 13440 \\ \hline
179 & $2^3 \cdot 3^2 \cdot 5 \cdot 89$ & 4 6 16 48 48 144 144 384 432 1152 3456 10368 15840 \\ \hline
181 & $2^3 \cdot 3^2 \cdot 5 \cdot 7 \cdot 13$ & 4 6 16 48 48 96 144 144 336 384 432 1152 3456 10368 16128 \\ \hline
191 & $2^7 \cdot 3 \cdot 5 \cdot 19$ & 4 6 16 24 48 48 96 192 384 720 768 1536 3072 12288 17280 \\ \hline
193 & $2^7 \cdot 3 \cdot 97$ & 4 6 16 24 48 96 192 384 768 1536 3072 12288 18816 \\ \hline
197 & $2^3 \cdot 3^2 \cdot 7^2 \cdot 11$ & 4 6 16 96 144 240 288 1920 4704 14112 17280 \\ \hline
199 & $2^4 \cdot 3^2 \cdot 5^2 \cdot 11$ & 4 6 16 24 48 144 144 240 576 1200 1920 3600 14400 17280 \\
\end{tabular}
\caption{Sizes of components under the action of Vieta moves, coordinate permutations, and double sign changes for the Cayley cubic $x^2 + y^2 + z^2 - xyz = 4$ mod $p$}
\label{table:cayley-components}
\end{table}

\begin{table}
\begin{tabular}{c|c|c}
$p$ & Factors of $p^2-1$ & {\rm Component sizes } \\ \hline
5 & $2^3 \cdot 3$ & 1 3 4 6 12 \\ \hline
7 & $2^4 \cdot 3$ & 1 3 4 6 12 24 \\ \hline
11 & $2^3 \cdot 3 \cdot 5$ & 1 3 4 6 12 12 36 48 \\ \hline
13 & $2^3 \cdot 3 \cdot 7$ & 1 3 4 6 12 24 48 72 \\ \hline
17 & $2^5 \cdot 3^2$ & 1 3 4 6 12 24 36 96 108 \\ \hline
19 & $2^3 \cdot 3^2 \cdot 5$ & 1 3 4 6 12 12 36 36 108 144 \\ \hline
23 & $2^4 \cdot 3 \cdot 11$ & 1 3 4 6 12 24 48 60 180 192 \\ \hline
29 & $2^3 \cdot 3 \cdot 5 \cdot 7$ & 1 3 4 6 12 12 24 36 72 96 288 288 \\ \hline
31 & $2^6 \cdot 3 \cdot 5$ & 1 3 4 6 12 12 24 36 96 96 288 384 \\ \hline
37 & $2^3 \cdot 3^2 \cdot 19$ & 1 3 4 6 12 36 48 108 180 432 540 \\ \hline
41 & $2^4 \cdot 3 \cdot 5 \cdot 7$ & 1 3 4 6 12 12 24 24 36 72 144 192 576 576 \\ \hline
43 & $2^3 \cdot 3 \cdot 7 \cdot 11$ & 1 3 4 6 12 24 60 72 180 192 576 720 \\ \hline
47 & $2^5 \cdot 3 \cdot 23$ & 1 3 4 6 12 24 48 96 192 264 768 792 \\ \hline
53 & $2^3 \cdot 3^3 \cdot 13$ & 1 3 4 6 12 36 84 108 252 324 972 1008 \\ \hline
59 & $2^3 \cdot 3 \cdot 5 \cdot 29$ & 1 3 4 6 12 12 36 48 96 144 288 420 1152 1260 \\ \hline
61 & $2^3 \cdot 3 \cdot 5 \cdot 31$ & 1 3 4 6 12 12 36 48 96 144 288 480 1152 1440 \\ \hline
67 & $2^3 \cdot 3 \cdot 11 \cdot 17$ & 1 3 4 6 12 60 144 180 432 480 1440 1728 \\ \hline
71 & $2^4 \cdot 3^2 \cdot 5 \cdot 7$ & 1 3 4 6 12 12 24 24 36 36 48 72 108 192 432 576 1728 1728 \\ \hline
73 & $2^4 \cdot 3^2 \cdot 37$ & 1 3 4 6 12 24 36 48 108 192 432 684 1728 2052 \\ \hline
79 & $2^5 \cdot 3 \cdot 5 \cdot 13$ & 1 3 4 6 12 12 24 36 84 96 144 252 576 672 2016 2304 \\ \hline
83 & $2^3 \cdot 3 \cdot 7 \cdot 41$ & 1 3 4 6 12 24 48 72 192 288 576 840 2304 2520 \\ \hline
89 & $2^4 \cdot 3^2 \cdot 5 \cdot 11$ & 1 3 4 6 12 12 24 36 36 60 96 108 180 288 720 864 2592 2880 \\ \hline
97 & $2^6 \cdot 3 \cdot 7^2$ & 1 3 4 6 12 24 24 48 72 96 192 384 768 1176 3072 3528 \\
\end{tabular}
\caption{Sizes of components under the action of Vieta moves, coordinate permutations for the Cayley cubic $x^2 + y^2 + z^2 - xyz = 4$ mod $p$}
\label{table:cayley-components-2-1}
\end{table}

\begin{table}
\begin{tabular}{c|c}
$p$ & {\rm Component sizes } \\ \hline
5 & 1 3 4 6 12 \\ \hline
7 & 1 3 4 6 12 24 \\ \hline
11 & 1 3 4 6 12 12 36 48 \\ \hline
13 & 1 3 4 6 12 24 48 72 \\ \hline
17 & 1 3 4 6 12 24 36 96 108 \\ \hline
19 & 1 3 4 6 12 12 36 36 108 144 \\ \hline
23 & 1 3 4 6 12 24 48 60 180 192 \\ \hline
29 & 1 3 4 6 12 12 24 36 72 96 288 288 \\ \hline
31 & 1 3 4 6 12 12 24 36 96 96 288 384 \\ \hline
37 & 1 3 4 6 12 36 48 108 180 432 540 \\ \hline
41 & 1 3 4 6 12 12 24 24 36 72 144 192 576 576 \\ \hline
43 & 1 3 4 6 12 24 60 72 180 192 576 720 \\ \hline
47 & 1 3 4 6 12 24 48 96 192 264 768 792 \\ \hline
53 & 1 3 4 6 12 36 84 108 252 324 972 1008 \\ \hline
59 & 1 3 4 6 12 12 36 48 96 144 288 420 1152 1260 \\ \hline
61 & 1 3 4 6 12 12 36 48 96 144 288 480 1152 1440 \\ \hline
67 & 1 3 4 6 12 60 144 180 432 480 1440 1728 \\ \hline
71 & 1 3 4 6 12 12 24 24 36 36 48 72 108 192 432 576 1728 1728 \\ \hline
73 & 1 3 4 6 12 24 36 48 108 192 432 684 1728 2052 \\ \hline
79 & 1 3 4 6 12 12 24 36 84 96 144 252 576 672 2016 2304 \\ \hline
83 & 1 3 4 6 12 24 48 72 192 288 576 840 2304 2520 \\ \hline
89 & 1 3 4 6 12 12 24 36 36 60 96 108 180 288 720 864 2592 2880 \\ \hline
97 & 1 3 4 6 12 24 24 48 72 96 192 384 768 1176 3072 3528 \\ \hline
101 & 1 3 4 6 12 12 36 144 144 300 432 900 1152 3456 3600 \\ \hline
103 & 1 3 4 6 12 24 84 144 252 432 1008 1152 3456 4032 \\ \hline
107 & 1 3 4 6 12 36 48 108 324 432 972 1404 3888 4212 \\ \hline
109 & 1 3 4 6 12 12 36 36 48 60 108 180 324 432 972 1440 3888 4320 \\ \hline
113 & 1 3 4 6 12 24 24 72 96 180 288 540 1152 1440 4320 4608 \\ \hline
127 & 1 3 4 6 12 24 24 36 72 96 108 192 384 576 1536 1728 5184 6144 \\ \hline
131 & 1 3 4 6 12 12 36 48 60 84 180 252 480 720 1440 2016 5760 6048 \\ \hline
137 & 1 3 4 6 12 24 144 264 432 792 1728 2112 6336 6912 \\ \hline
139 & 1 3 4 6 12 12 24 36 72 144 264 288 576 792 1728 2112 6336 6912 \\ \hline
149 & 1 3 4 6 12 12 36 96 288 300 684 900 2052 2400 7200 8208 \\ \hline
151 & 1 3 4 6 12 12 24 36 96 180 288 300 540 900 2160 2400 7200 8640 \\ \hline
157 & 1 3 4 6 12 48 84 252 672 1008 2016 3120 8064 9360 \\ \hline
163 & 1 3 4 6 12 36 108 324 840 972 2520 2916 8748 10080 \\ \hline
167 & 1 3 4 6 12 24 24 48 72 192 192 288 576 1152 2304 3444 9216 10332 \\ \hline
173 & 1 3 4 6 12 420 924 1260 2772 3360 10080 11088 \\ \hline
179 & 1 3 4 6 12 12 36 36 48 96 108 144 288 432 864 1152 2592 3960 10368 11880 \\
\end{tabular}
\caption{Sizes of components under the action of Vieta moves, coordinate permutations for the Cayley cubic $x^2 + y^2 + z^2 - xyz = 4$ mod $p$}
\label{table:cayley-components-2-2}
\end{table}

We note that, modulo any prime $p \geq 5$, the Cayley cubic has components of size 4, 6, and 16. These come from particular solutions over $\Z$. The component of size 4  consists of $(2,2,2)$ and its orbit under double sign changes, namely $(-2,-2,2), (-2, 2, -2)$, and $(2, -2, -2)$. The orbit of size 6 consists of permutations of $(0,0,\pm2)$, which is a special case of the size 6 component that arises from $(0,0,\pm \sqrt{k})$ whenever $k$ is a square modulo $p$. The orbit of size 16 consists of $(1,1,2)$, its Vieta image $(1,1,-1)$, and their orbits under permutations and double sign changes. Using only Markoff moves and permutations, without sign changes, one would have $(1+3)+6+(12+4)$ instead of $4+6+16$: Namely $(2,2,2)$ in its own orbit, an orbit of size 3 containing $(2,-2,-2)$, the orbit of size 6 containing $(2,0,0)$, an orbit of size 12 containing $(1,1,2)$, and another orbit of size 4 containing $(-1,-1,-1)$.

For many primes $p$, there is a component of size 24, and this also has a simple explanation. If 2 is a square modulo $p$, then among the solutions to equation (\ref{eqn:cayley}) are $(\sqrt{2},\sqrt{2},0)$ and its Vieta image $(\sqrt{2},\sqrt{2},2)$. Permutations and double sign changes of these then yield a component of size 24. It consists of the vectors $(\varepsilon_1 \sqrt{2}, \varepsilon_2 \sqrt{2},0)$, $(\varepsilon \sqrt{2}, \varepsilon \sqrt{2},2)$, $(\varepsilon \sqrt{2}, -\varepsilon \sqrt{2}, -2)$ and their permutations, where $\varepsilon, \varepsilon_1, \varepsilon_2 = \pm 1$ are signs. These can also be reached from one another using Markoff moves instead of sign changes. By the supplement to the law of quadratic reciprocity, 2 is a square if and only if $p^2 -1$ is divisible by 16:
\begin{equation} \label{eqn:supplement}
\legendre{2}{p} = (-1)^{(p^2-1)/8}.
\end{equation}
This explains why, in Table~\ref{table:cayley-components}, the primes 7, 17, 23, 31, 41, 47, 71, 73, 79, 89, and 97 are precisely the ones with a component of size 24. It is also a clue that the other component sizes might be explained most directly in terms of $p^2 - 1$ and its factors.

A special feature of the equation $x^2+y^2+z^2 = xyz$ is that when we change variables to 
\begin{equation} \label{eqn:2cosh}
x = \xi + \xi^{-1}, \quad y = \eta + \eta^{-1}
\end{equation}
the solutions for $z$ are then
\begin{equation} \label{eqn:zsol}
\xi \eta + \frac{1}{\xi \eta}, \quad \xi \eta^{-1} + \frac{1}{\xi \eta^{-1} }.
\end{equation}
For $x \in \F_p$, there is a solution $\xi \in \F_p^{\times}$ when $x^2-4$ is a square mod $p$. Otherwise, $\xi$ must be taken from a quadratic extension $\F_{p^2}$. Thus we let $g$ be a generator of $\F_{p^2}^{\times}$ and write
\begin{equation} \label{eqn:g2cosh}
x = g^u + g^{-u}, \quad y = g^v + g^{-v}
\end{equation}
where the exponents are taken modulo $p^2-1$. 
The solutions for $z$ are then 
\begin{equation} \label{eqn:gzsol}
g^{u+v} + g^{-u-v}, \quad g^{u-v} + g^{-u+v}.
\end{equation}
Note that $-u$ and $u$ define the same $x$, and likewise $v$ is equivalent to $-v$. 
Hence $(u,v,u-v)$ is equivalent to $(u,-v, u + (-v))$, so that all solutions can be parametrized in the form $(u,v,u+v)$ with the third coordinate equal to the sum of the others.
If we had chosen a different generator, say $g^w$ instead of $g$, the exponents $u$ and $v$ would simply be multiplied by a unit $w$ modulo $p^2-1$. 
We are interested in the ``real" solutions to \ref{eqn:cayley}, that is to say those over $\F_p$ rather than $\F_{p^2}$. To have $x = g^u + g^{-u}$ lie in $\F_p$, it is necessary and sufficient that it be fixed by the Galois involution $x \mapsto x^p$. This holds if and only if $pu \equiv \pm u \bmod p^2-1$. Thus $u$ must be a multiple of $p+1$ or $p-1$. Likewise, the second coordinate $v$ must be a multiple of $p+1$ or $p-1$, perhaps not the same one as $u$. If both $u$ and $v$ are multiples of the same $p \pm 1$, then $p \pm 1$ also divides the sum $u+v$ and so the third coordinate is also real.

\begin{proposition} \label{prop:matrices}
The Markoff moves, as well as coordinate permutations, act linearly on the coordinates $\begin{bmatrix} u \\ v \end{bmatrix}$. Explicitly, their matrices are given by
\begin{equation}
[m_1] = \begin{bmatrix} 1 & 2 \\ 0 & -1 \end{bmatrix}, \quad [m_2] = \begin{bmatrix} 1 & 0 \\ -2 & - 1 \end{bmatrix}, \quad [m_3] = \begin{bmatrix} 1 & 0 \\ 0 & -1 \end{bmatrix}.
\label{eqn:markoff-matrices}
\end{equation}
and
\begin{equation}
[\tau_{12}] = \begin{bmatrix} 0 & 1 \\ 1 & 0 \end{bmatrix}, [\tau_{23}] = \begin{bmatrix} -1 & 0 \\ 1 & 1 \end{bmatrix}, 
[\tau_{13} ] = \begin{bmatrix} 1 & 1 \\ 0 & -1 \end{bmatrix}
\label{eqn:transpo-matrices}
\end{equation}
where $\tau_{ij}$ is the transposition exchanging $i$ and $j$.
These matrices generate $\GL_2(\Z)$ or, modulo $p^2-1$, the subgroup of matrices with determinant $\pm 1$.
\end{proposition}
Note that these matrices are better interpreted in $\PGL_2$ than $\GL_2$ because the exponents $u$ and $v$ are only defined up to sign. One must change the sign of the entire vector because changing the sign of only one of $u, v$ will not keep the third coordinate $u+v$ equal to the sum of the others. 
\begin{proof}
First note that $m_3$ exchanges $(u,v,u+v)$ with $(u,v,u-v)$, or equivalently with $(u,-v, u -v)$. 
We use the latter form to keep the third coordinate equal to the sum of the first two.
The transposition $\tau_{12}$ sends $(u,v,u+v)$ to $(v,u,u+v)$.
The transposition $\tau_{23}$ sends $(u,v,u+v)$ to $(u,u+v,v)$, or equivalently $(-u,u+v,-u + (u+v))$.
The transposition $\tau_{13}$ sends $(u,v,u+v)$ to $(u+v,v,u)$ or equivalently $(u+v,-v,u)$.
In matrix form acting on $\begin{bmatrix} u \\ v \end{bmatrix}$, these operations correspond to
\[
[m_3] = \begin{bmatrix} 1 & 0 \\ 0 & -1 \end{bmatrix},
[\tau_{12}] = \begin{bmatrix} 0 & 1 \\ 1 & 0 \end{bmatrix},
[\tau_{23}] = \begin{bmatrix} -1 & 0 \\ 1 & 1 \end{bmatrix},
[\tau_{13} ] = \begin{bmatrix} 1 & 1 \\ 0 & -1 \end{bmatrix}.
\]
Using the relations $m_2 = \tau_{23} m_3 \tau_{23}$ and $m_1 = \tau_{13} m_3 \tau_{13}$, we then find
\[
[m_2] = \begin{bmatrix} 1 & 0 \\ -2 & - 1 \end{bmatrix}, \quad [m_1] = \begin{bmatrix} 1 & 2 \\ 0 & -1 \end{bmatrix}.
\]
To determine what group these matrices generate, note that multiplying by $[m_3] = \begin{bmatrix} 1 & 0 \\ 0 & -1 \end{bmatrix}$ changes the sign of the second row or column:
\[
\begin{bmatrix} a & b \\ c & d \end{bmatrix} [m_3] = \begin{bmatrix} a & -b \\ c & -d \end{bmatrix}, \quad [m_3]\begin{bmatrix} a & b \\ c & d \end{bmatrix} = \begin{bmatrix} a & b \\ -c & -d \end{bmatrix}.
\]
Combining this with $\tau_{23} = \begin{bmatrix} 0 & 1 \\ 1 & 0 \end{bmatrix}$, which exchanges two rows or columns, we may also change the sign of the first row or column.
This is enough to obtain the standard generators for $\SL_2(\Z)$:
\begin{align*}
T &= \begin{bmatrix} 1 & 1 \\ 0 & 1 \end{bmatrix} = [m_3] [\tau_{13}] \\
S &= \begin{bmatrix} 0 & 1 \\ -1 & 0 \end{bmatrix} = [m_3][\tau_{12}]
\end{align*}
One also has $S = -[\tau_{12}][m_3]$. Hence the group generated by the matrices (\ref{eqn:markoff-matrices}) and (\ref{eqn:transpo-matrices}) contains $\SL_2(\Z)$. Multiplying by any matrix of determinant $-1$, for instance $\tau_{23}$, we obtain the other coset of $\SL_2(\Z)$ in $\GL_2(\Z)$. Hence these matrices generate $\GL_2(\Z)$.
\end{proof}

To obtain simpler graphs, we have previously used $D_1 = \tau_{23} \circ m_1 , D_2 = \tau_{23} \circ m_2, D_3= \tau_{23} \circ m_3$, which do not generate all of $\GL_2(\Z)$. But for Table~\ref{table:cayley-components}, we have used the full symmetry of all the Markoff moves, all the transpositions, and also double sign changes.
The double sign changes do not act linearly on the exponents $(u,v)$. Instead, since
\begin{equation}
-1 = g^{ \frac{p^2-1}{2} } = g^{- \frac{p^2-1}{2} }
\end{equation}
their effect is to translate one or both of $u, v$ by $(p^2-1)/2$.
Note that, modulo $p^2-1$, the exponent for the third coordinate remains equal to $u+v$: It is translated by $(p^2-1)/2$ if only one of $u,v$ is, or by $(p^2-1)/2 + (p^2-1)/2 = 0$ if both are.
We will first determine the orbits under the linear action, and then incorporate these three translations. 
The linear action is dictated by matrix arithmetic modulo $p^2-1$, which can be understood via the Chinese remainder theorem and the corresponding action modulo prime powers. This is the underlying reason that the factors of $p^2-1$ play such an important role in the structure of the Cayley cubic.

\subsection{Proof of Theorem~\ref{thm:main}: Number of orbits}
Consider the action of matrices with determinant $\pm 1$ on $\Z/q^e \times \Z/q^e$, where $q^e$ is a prime power.
Given a vector $\begin{bmatrix} a \\ c \end{bmatrix}$ where at least one of $a, c$ is invertible, either
\[
\begin{bmatrix} a & 0 \\ c & a^{-1} \end{bmatrix} \ \text{or} \ \begin{bmatrix} a & -c^{-1} \\ c & 0 \end{bmatrix}
\]
will have determinant 1 and send $\begin{bmatrix} 1 \\ 0 \end{bmatrix}$ to $\begin{bmatrix} a \\ c \end{bmatrix}$.
If neither $a$ nor $c$ is invertible modulo $q^e$, then they must be divisible by $q$. Let $q^f$ be the largest power of $q$ dividing both of them. Note that $A(q^f w) = q^f Aw$, so that every vector in the orbit of $w$ also has both coordinates divisible by $q^f$. Conversely, since $q^f$ is the largest power of $q$ dividing both, either $w_1/q^f$ or $w_2/f$ is a unit. Thus there is a matrix of determinant 1 taking $\begin{bmatrix} 1 \\ 0 \end{bmatrix}$ to $w/q^f$, which shows that $w$ is in the same orbit as $\begin{bmatrix} q^f \\ 0 \end{bmatrix}$.
It follows that $\SL_2(q^e)$ has $e+1$ orbits on $\Z/q^e \times \Z/q^e$ and a list of representatives is
\[
\begin{bmatrix} 1 \\ 0 \end{bmatrix}, \begin{bmatrix} q \\ 0 \end{bmatrix}, \ldots, \begin{bmatrix} q^{e-1} \\ 0 \end{bmatrix}, \begin{bmatrix} q^e \\ 0 \end{bmatrix} = \begin{bmatrix} 0 \\ 0 \end{bmatrix}.
\]
The orbits are the same under the group of matrices of determinant $\pm 1$ or even the full group $\GL_2(\Z/q^e)$.
Modulo a composite $N$, two vectors are in the same orbit if and only if their images modulo $q^e$ are in the same orbit for each prime power factor of $N$. 
The orbits for the action of $\SL_2(\Z/N)$ on $\Z/N \times \Z/N$ can be found by the Chinese remainder theorem, and likewise for $\GL_2(\Z/N)$ or the subgroup of matrices with determinant $\pm 1 \bmod N$. 
The orbits are parameterized by all choices of $\{ f(q) \}_{q|N}$, where $0 \leq f(q) \leq e_q$ specifies the highest power of $q$ that divides the coordinates of vectors in a given orbit.
Equivalently, we may think of the parameter $f$ as a divisor of $N$, namely $t = \prod q^{f(q)}$, and then the corresponding orbit simply consists of vectors both of whose coordinates are divisible by $t$.
From either perspective, the number of orbits is therefore
\[
\prod_{q |N} (e_q+1)
\]
where $q$ ranges over all prime divisors of $N$ and $e_q$ is the highest power of $q$ dividing $N$. 

With $N = p^2-1$, all of these orbits are candidates as orbits for the action of permutations and Markoff moves on the Cayley cubic. 
However, if the coordinates are not divisible by $p+1$ or $p-1$, one obtains solutions over the extension $\F_{p^2}$ rather than $\F_p$. We must discard these orbits.
We must also identify $(u,v)$ and $(-u,-v)$ because they define the same solution $(x,y,z)$, but this does not change the number of orbits because each orbit is already closed under negation.

Note that $p-1$ and $p+1$ have no common factor except $2$. If $p \equiv 1 \bmod 4$, then $p-1$ is ``highly" divisible by 2 while $p+1$ is only once divisibe by 2. If $p \equiv -1 \bmod 4$, then it is $p+1$ that contains most of the factors of 2. To avoid considering these cases separately, let $p+\varepsilon$ be divisible simply by 2 and $p- \varepsilon$ by the remaining factors of 2. Here, the sign is
\begin{equation}
\varepsilon = (-1)^{\frac{p-1}{2}} = \legendre{-1}{p}.
\end{equation}
Let $Q^+$ and $Q^-$ be the sets of odd primes dividing $p+\varepsilon$ or $p - \varepsilon$ respectively. These are disjoint. Thus
\begin{equation}
p+\varepsilon = 2 \prod_{q \in Q^+} q^{e_q}, \quad p - \varepsilon = 2^{e_2 -1} \prod_{q \in Q^{-}} q^{e_q}.
\end{equation}
The ``real" orbits are the ones with either
\begin{itemize}
\item
$f(2) \geq 1$ and $f(q) = e_q$ for all $q \in Q^+$, or 
\item
$f(2) \geq e-1$ and $f(q) = e_q$ for all $q \in Q^-$, 
\end{itemize}
or both.
In the first case, $f(q)$ assumes any of $e_2$ values $1, 2, \ldots, e_2$ for $q=2$, must equal $e_q$ for $q \in Q^+$, and could be any of $0, 1, \ldots, e_q$ for $q \in Q^-$. In the second case, $f(2)$ takes only two values $e_2 - 1$ or $e_2$, $f(q)$ must equal $e_q$ for $q \in Q^-$, and could be any of $0, 1, \ldots, e_q$ for $q \in Q^+$. In case of overlap, both coordinates are divisible by $\text{lcm}(p+\varepsilon,p-\varepsilon) = (p^2-1)/2$. This only happens for two orbits, namely $f(2)$ may be $e_2$ or $e_2-1$ but $f(q)$ must equal $e_q$ for all odd $q$.  
We subtract 2 to compensate for double-counting these two orbits. The total is
\begin{equation}
e_2  \prod_{q \in Q^-} (e_q+1) + 2 \prod_{q \in Q^+} (e_q+1)  - 2.
\end{equation}
This is the formula stated in Theorem~\ref{thm:main}. \qed

\subsection{Including double sign changes}
Now we incorporate the further symmetry of the Cayley cubic under sign changes of the form $(x,y,z) \mapsto (\varepsilon_1 x, \varepsilon_2 y, \varepsilon_3 z)$ with $\varepsilon_1 \varepsilon_2 \varepsilon_3 = 1$.
The double sign change $\sigma_{12} = (x,y,z) \mapsto (-x,-y,z)$ acts on the exponents $(u,v)$ by
\begin{equation}
\sigma_{12} : (u,v) \mapsto \left( u + \frac{p^2-1}{2}, v + \frac{p^2-1}{2} \right)
\end{equation}
because $-1 = g^{(p^2-1)/2}$. Note that, working modulo $p^2-1$, the exponent for $z = g^{u+v} + g^{-u-v}$ remains the sum of the exponents for $x$ and $y$. The other sign changes are conjugate to this one by transpositions:
\begin{align*}
\sigma_{13} &= \tau_{23} \sigma_{12} \tau_{23} \\
\sigma_{23} &= \tau_{13} \sigma_{12} \tau_{13}.
\end{align*}
Therefore it is enough to determine how the Markoff+permutation orbits above merge under the action of $\sigma_{12}$. 
For the odd primes $q$ dividing $p^2-1$, note that $(p^2-1)/2$ remains equally divisible by $q$, so $\sigma_{12}$ acts trivially modulo $q^e$.
For $q=2$, note that $(p^2-1)/2$ is only divisible by $2^{e-1}$ instead of $2^e$. 
Thus the orbits where $f(2) < e_2-1$ are not affected, but an orbit with $f(2) = e_2-1$ merges with the orbit having $f(2) = e_2$ and the same value $f(q)$ for odd $q$.
The factor of 2 must be removed in the product over $Q^+$, because the orbits divisible by $p - \varepsilon$ merge pairwise.
The orbits divisible by $p + \varepsilon$ are not affected, unless $f(2) = e_2-1$ or $e_2$. Effectively, there is one less choice for $f(2)$ so the factor $e_2$ is replaced by $e_2-1$. 
The two ``overlap" orbits divisible by $(p^2-1)/2$ have been removed twice in this process, so we must add 1 to compensate. 
We must therefore subtract 1 instead of 2 compared to the formula above, because now only one orbit is double-counted.
The total is then
\begin{equation}
(e_2 -1 ) \prod_{q \in Q^-} (e_q+1) +  \prod_{q \in Q^+} (e_q+1) - 1
\end{equation}
as stated in Corollary~\ref{cor:signs}, (\ref{eqn:num-orb-sign}).
\qed

\subsection{Proof of Theorem~\ref{thm:main}: Sizes of orbits}
Let us determine the size of the orbit of $\begin{bmatrix} t \\ 0 \end{bmatrix}$, where $t = \prod q^{f(q)}$ is a divisor of $N = p^2-1$. For any group action, we have the orbit-stabilizer formula
\begin{equation}
\# \Orb(t) = \frac{\# G}{\# \Stab(t)}.
\end{equation}
In the present case, the stabilizer consists of matrices sending $\begin{bmatrix} t \\ 0 \end{bmatrix}$ to itself or alternatively to $\begin{bmatrix} -t \\ 0 \end{bmatrix}$, since these define the same solution to (\ref{eqn:cayley}).
The initial group consists of matrices with determinant $\pm 1 \bmod N$, but the orbit structure is the same for $\GL_2(\Z/N)$ and we will make the replacement $G = \GL_2$ to avoid enforcing the determinant condition when we determine $\Stab(t)$ and to simplify the expression for $\# G$. 
As a base case, $|\GL_2(\Z/q\Z) | = (q^2-1)(q^2 - q)$ because there are $q^2 - 1$ non-zero choices for the first column, and then $q^2- q$ vectors not equal to a multiple of the first column.
To pass to higher powers of $q$, we use a version of Hensel's Lemma to lift the matrices from $\GL_2(\Z/q\Z)$. Write the matrix whose invertibility is to be determined as $A + qA' + q^2 A'' + \ldots + q^{e-1} A^{(e-1)}$, where each of the matrices $A, A', \ldots$ has entries in $\Z/q\Z$. The condition for another such matrix $B + qB' + q^2 B'' + \ldots$ to be its inverse is that
\begin{align*}
I &= (A + qA' + \ldots)(B + qB' + \ldots) \\
&= AB + q(A'B + AB') + q^2(A' B' + A'' B + AB'') + \ldots
\end{align*}
Thus we must first of all have $AB = I$, so that $A$ is in $\GL_2(\F_q)$. Then we must have $A'B + AB' = 0$, which can be arranged for any choice of $B'$ by taking $A' = -AB'B^{-1} = -AB'A$. Thus we have $q^4$ choices at this stage. Then we must have $A'B' + A''B + AB'' = 0$, which determines $A'' = -(AB'' + A'B')B^{-1}$ given any of $q^4$ choices for $B''$.
Continuing in this way, we find that
\[
|\GL_2(\Z/q^e\Z) | = |\GL_2(\F_q)| q^{4(e-1)}.
\]
For $N = \prod_{q^e | N} q^e$, it follows that
\begin{equation}
| \GL_2(\Z/N\Z) | = \prod_{q^e | N} (q^2-1)(q^2-q)q^{4(e-1)} .
\end{equation}

To determine the stabilizer, we suppose that
\begin{equation}
A\begin{bmatrix} t \\ 0 \end{bmatrix} \equiv \pm \begin{bmatrix}  t \\ 0 \end{bmatrix} \bmod N.
\end{equation}
The same congruence holds modulo each prime power $q^e$ (with the same choice of $\pm$), or equivalently
\begin{equation}
A \begin{bmatrix} 1 \\ 0 \end{bmatrix} \equiv \begin{bmatrix}  \pm 1 \\ 0 \end{bmatrix} \bmod q^{e-f}
\end{equation}
where $f$ is the highest power of $q$ dividing $t$. If $e=f$, then $t = 0$ and there is no constraint on $A$. If $f(q) < e_q$, then the first column of $A$ is constrained. To ensure invertibiility, the second diagonal entry must be a unit modulo $q^{e-f}$. When we lift $A$ to $\Z/q^e$, it must take the form
\begin{equation}
A \in \begin{bmatrix} \pm 1 & \Z/q^{e-f} \\ 0 & {\Z/q^{e-f} }^{\times} \end{bmatrix} + q^{e-f} M_{2\times 2}(\Z/q) + \ldots + q^{e-1} M_{2\times 2}(\Z/q).
\end{equation}
There is only one choice for the first column, since the sign $\pm$ is fixed and common to all factors $q^e$.
There are $\phi(q^{e-f})$ choices for the second diagonal entry, $q^{e-f}$ choices for the other entry in the second column, and $q^{4f}$ ways to lift. It follows that, in the action on $\Z/N \times \Z/N$, the stabilizer has size
\begin{equation}
\# \Stab(t) = \prod_{q: f(q) = e_q} |\GL_2(\Z/q^e)| \times \prod_{f(q) < e_q} \phi(q^{e-f}) q^{e-f} q^{4f} .
\end{equation}
In the action on the Cayley cubic, the stabilizer is usually twice this size because $(-t,0)$ represents the same solution as $(t,0)$. The exceptional cases are $t = 0$ and $t = (p^2-1)/2$, since then $t = -t$.
Note that
\[
\frac{|\GL_2(\Z/q^e)|}{\phi(q^{e-f}) q^{e-f} q^{4f} } = \frac{(q^2-1)(q^2-q) q^{4(e-1)} }{2q^{e-f-1}(q-1)q^{e-3f}} = (q^2-1) q^{2e-2f-2}.
\]
By the orbit-stabilizer formula, the size of the corresponding orbit is
\begin{equation}
\# \Orb(t) = \frac{1}{2} \prod_{f(q) < e_q} (q^2-1) q^{2e-2f -2}
\end{equation}
except without the factor $1/2$ if $t=0$ or $t = (p^2-1)/2$. This completes the proof of Theorem~\ref{thm:main}. \qed

\subsection{Examples: Sizes of orbits}
Note that $p^2-1$ is divisible by 8 and by 3 for any odd $p$, and hence modulo any $p$ there are divisors
\begin{equation}
t = 0, \frac{p^2-1}{2}, \frac{p^2-1}{2}, \frac{p^2-1}{3}, \frac{p^2-1}{4}, \frac{p^2-1}{3}, \frac{p^2-1}{6}
\end{equation}
as well as
\begin{equation}
t = \frac{p^2-1}{8}, \frac{p^2-1}{12}, \frac{p^2-1}{24}.
\end{equation}
We have listed these separately because the divisors in the first list are automatically divisible by $p+1$ or $p-1$, while those in the second list may or may not be.
We start with 
\begin{align*}
g^0 + g^{-0} &= 1 + 1 = 2 \\
g^{\frac{p^2-1}{2}} + g^{-\frac{p^2-1}{2}} &= -1 -1 = -2.
\end{align*} 
We solve for the others using ``bisection", that is, substituting previously known values into the relation
\begin{equation}
\big( g^{u/2} + g^{-u/2} \big)^2 = g^u + g^{-u} + 2.
\end{equation}
For example, $x = g^{(p^2-1)/4} + g^{-(p^2-1)/4}$ solves $x^2 = 0$ so it must be that $x = 0$.
We have $(p^2-1)/3 \equiv -2(p^2-1)/3 \bmod p^2-1$, so $y = g^{(p^2-1)/3} + g^{-(p^2-1)/3}$ solves $y^2 = y+2$. Therefore $y=-1$, since 0 is already spoken for. Then we simply multiply by $-1 = g^{(p^2-1)/2}$ to find
\[
g^{\frac{p^2-1}{6}} + g^{-\frac{p^2-1}{6}} = 1
\]
or alternatively solve the equation $z^2 = -1 + 2$ using the previous value for $-1$. 
``Bisecting" these values, we find that
\begin{align*}
g^{\frac{p^2-1}{8}} + g^{-\frac{p^2-1}{8}} &= \sqrt{2} \\
g^{\frac{p^2-1}{12}} + g^{-\frac{p^2-1}{12}} &= \sqrt{3} \\
g^{\frac{p^2-1}{24}} + g^{-\frac{p^2-1}{24}} &= \sqrt{2+\sqrt{3}}
\end{align*}
which may or may not lie in $\F_p$. In any case, we can determine the size of the corresponding orbit in $\Z/N \times \Z/N$ and, if the necessary coordinates lie in $\F_p$, the Cayley cubic will have an orbit of this size.

The size of the orbit corresponding to a divisor $t = \prod q^{f(q)}$ is
\begin{equation}
\frac{1}{2} \prod_{q: f(q) < e_q} (q^2-1) q^{2(e-f-1)}
\end{equation}
or twice that in case $t = 0$ or $t = \frac{p^2-1}{2}$. 
The easiest case is the orbit of $(u,v)=(0,0)$, which obviously has size 1. This is the case $t=0$ and our formula also gives 1, because the product is empty and the factor $1/2$ is omitted. This is the orbit of $(x,y,z)=(2,2,2)$ in the original coordinates. For $t = (p^2-1)/2$, i.e. the orbit of $(-2,2,-2)$, we again omit the factor $1/2$ and find that the orbit has size $(2^2 - 1) \cdot 2^0 = 3$.
For $t = (p^2-1)/4$, i.e. the orbit of $(0,2,0)$, we have $f(2)=e_2 - 2$ so the orbit size is
\[
\frac{1}{2} (2^2 -1) \cdot 2^{2(2-1)} = 6.
\]
For $t = (p^2-1)/3$, i.e. the orbit of $(-1,2,-1)$, we have $f(3) = e_3-1$ and $f(q)=e_q$ otherwise, so the size of the orbit is
\[
\frac{1}{2} (3^2-1) \cdot 3^0 = 4.
\]
For $t=(p^2-1)/6$, i.e. the orbit of $(1,2,1)$, we have $f(2) = e_2-1$ as well as $f(3) = e_3-1$ so the orbit size is
\[
\frac{1}{2} (2^2-1) \cdot 2^0 \times (3^2-1) \cdot 3^0 = 12.
\]
This is another way to explain the orbits of size 1, 3, 4, 6, 12 which are present modulo any prime.
Recall that double sign changes merge the orbit corresponding to $t$ with the orbit corresponding to $2t$ whenever $f(2) = e_2-1$. Thus the orbits of $(p^2-1)/2$ and $0$ merge, as do the orbits of $(p^2-1)/6$ and $(p^2-1)/3$.
This is why the sizes 4, 6, 16 appear in Table~\ref{table:cayley-components}.

If $p^2-1$ is divisible by 16, then we also have the orbit of $(\sqrt{2},2,\sqrt{2})$ with $t=(p^2-1)/8$. Because $f(2) = e_2 - 3$, the size of this orbit is
\[
\frac{1}{2} (2^2-1) 2^{2(3-1)} = 24.
\]
If 3 is a square mod $p$ and we take $t= (p^2-1)/12$, then we have an orbit with $f(2) = e_2-2$ and $f(3)=e_3-1$ and hence of size
\[
\frac{1}{2} (2^2-1) \cdot 2^{2(2-1)} \times \frac{3^2-1}{2} \cdot 3^0 = 48.
\]
This occurs when $p \equiv \pm 1 \bmod 12$, by quadratic reciprocity.
If both 2 and 3 are squares, then for $t=(p^2-1)/24$ we have $f(2)=e_2-3$ and $f(3)=e_3-1$, which gives an orbit of size
\[
\frac{1}{2} (2^2-1) \cdot 2^{2(3-1)} \times (3^2-1) = 3 \cdot 2^6 = 192.
\]
This component first occurs when $p=23$.
None of these components merge under double sign changes, because $f(2) < e_2 - 1$. 

\subsection{Examples: Number of orbits}
First, consider the case without sign changes. When $p=29$, we have $p-1 = 2^2 \cdot 7$ and $p+1 = 2 \cdot 3 \cdot 5$, so $Q^+$ consists of 3 and 5 while $Q^-$ consists of 7. The exponent $e_2$ is 3. The formula (\ref{eqn:num-orb}) gives
\[
3 \times 2 + 2 \times (2 \times 2)  - 2 = 12.
\]
For example, when $p = 71$, we have $p-1 = 2 \cdot 5 \cdot 7$ and $p+1 = 2^3 \cdot 3^2$, so $Q^+ = \{ 5, 7 \}$ and $Q^- = \{ 3 \}$. The formula gives
\[
e_2  \prod_{q \in Q^-} (e_q+1) + 2 \prod_{q \in Q^+} (e_q+1)  - 2 = 4 \times 3 + 2 \times (2 \times 2) - 2 = 18
\]
and indeed there are 18 orbits (of respective sizes 1 3 4 6 12 12 24 24 36 36 48 72 108 192 432 576 1728 1728).

For a first example including sign changes, take $p = 5$. Then $Q^-$ is empty, $Q^+ = \{ 3 \}$, and $e_2 = 3$, so the Cayley cubic splits into $2+2-1 = 3$ orbits. 
When $p=199$, we have $p-1 = 2^2 \cdot 3^2 \cdot 11$ and $p+1 = 2^3 \cdot 5^2$, so $Q^0 = \{ 5 \}$, $Q^+ = \{3, 11 \}$, and $e_2 = 4$. We have $e_3 = 2 = e_5$ and $e_{11}=1$, so the number of orbits is $3 \times 3 + 3\times 2 - 1 = 14$. This explains the number of orbits in Table~\ref{table:cayley-components}.

As a final example, suppose $p = 2l+1$ is a Sophie Germain prime. Then $Q^+$ consists only of the prime $l$ and $Q^-$ contains the prime factors of $l+1$. The number of orbits (including sign changes) is
\begin{equation}
\text{ord}_2(l+1) \prod_{q | (l+1)} \big( \text{ord}_q(l+1) + 1 \big).
\end{equation}

\subsection{Finite orbits in characteristic 0}
The finite orbits over $\C$ are determined by roots of unity. Whenever the finite field $\F_p$ contains a particular root of unity, the corresponding orbit will appear in the Cayley cubic mod $p$.
Suppose $(x,y,z)$ belongs to a finite orbit of the Cayley cubic over $\C$. Then some power of the element $\tau_{23} \circ m_3$ must take $(x,y,z)$ to itself. We have
\[
\tau_{23} \circ m_3: (x,y,z) \mapsto (x,y,xy-z) \mapsto (x,xy-z,y).
\]
Thus the latter two coordinates are transformed by the matrix $\begin{bmatrix} x & -1 \\ 1 & 0 \end{bmatrix}$, which must have finite order if we are to return to $(x,y,z)$ after finitely many steps. Thus its eigenvalues must be roots of unity. The trace is $x$ and the determinant is $1$, so the eigenvalues are $\xi, \xi^{-1}$ where $x = \xi+\xi^{-1}$. To have $\xi^n = 1$, we must have $x = e^{2\pi i k/n}$ for some integer $k$. Then
\[
x = \xi + \xi^{-1} = 2 \cos \frac{2\pi k}{n}.
\]
A similar conclusion for $y$ follows by considering $\tau_{23} \circ m_2$, which acts on $(y,z)$ by $\begin{bmatrix} 0 & 1 \\ -1 & x \end{bmatrix}$. Then using (\ref{eqn:cayley}), we deduce from $x = \xi + \xi^{-1}$ and $y = \eta + \eta^{-1}$ that
\[
z = \xi \eta + \frac{1}{\xi \eta}, \ \text{or} \ z = \xi^{-1}\eta + \frac{1}{\xi^{-1}\eta}
\]
Thus for $(x,y,z)$ to be part of a finite orbit, it must be of the form 
\[
2( \cos{\alpha}, \cos{\beta}, \cos(\alpha \pm \beta) )
\]
where $\alpha, \beta$ are rational multiples of $\pi$. Conversely, applying Markoff moves and permutations to such a point will not increase the denominators of the angles $\alpha, \beta$, so its orbit will be finite. Dubrovin and Mazzocco do a similar calculation in the context of braid groups in \cite{DubMaz}, Lemma 1.12.

\subsection{Comparison with other levels}
All level sets $x^2+y^2+z^2 = xyz + k$ have an interpretation by which the same matrices from Proposition~\ref{prop:matrices} act. This is given by \emph{Fricke's trace identity}. For $2 \times 2$ matrices of determinant 1,
\begin{equation} \label{eqn:fricke}
\tr(A)^2 + \tr(B)^2 + \tr(AB)^2 = \tr(A)\tr(B)\tr(AB) + \tr(ABA^{-1}B^{-1}) + 2.
\end{equation}
Thus if $\tr(ABA^{-1}B^{-1}) = k-2$, the vector of traces $(\tr{A}, \tr{B}, \tr{AB})$ solves the Markoff equation at level $k$. If $AB = BA$, then $\tr(ABA^{-1}B^{-1}) = 2$ and we have a point on the Cayley cubic.
For a matrix of determinant 1, the eigenvalues form a pair $\xi, \xi^{-1}$ inverse to each other, so the trace is
\begin{equation}
\tr(A) = \xi + \xi^{-1}
\end{equation}
and this is the same change of variable from the beginning of this section.
If instead $AB = -BA$, then $\tr(ABA^{-1}B^{-1}) = -2$ and we obtain points on the original Markoff surface at level $k=0$. If $Av = \xi v$, then anticommuting with $B$ forces $A(Bv) = -\xi(Bv)$ so that $-\xi$ is also an eigenvalue. To have $\det(A) = 1$, this implies that $\xi^2 = -1$. Likewise, the eigenvalues of $B$ must be $\pm \sqrt{-1}$. If $p \equiv 1 \bmod 4$, then the ground field contains a $\sqrt{-1}$ and solutions of this form are helpful in constructing the giant component of Bourgain-Gamburd-Sarnak \cite{BGS}.

Both sides of (\ref{eqn:fricke}) are polynomials in the eight entries of the two matrices, since the determinant being 1 allows one to skip the division by $\det(A), \det(B)$ in computing $A^{-1}, B^{-1}$. In principle, one can manually verify that they coincide. A more elegant proof is made possible by the Cayley-Hamilton theorem, the cyclic property $\tr(XY) = \tr(YX)$, and the fact that $\tr(X) = \tr(X^{-1})$ for $X \in \SL_2$. See \cite{A}, Proposition 4.3, p. 65. The argument is related to why the matrices in Proposition~\ref{prop:matrices} give the action of Markoff moves and permutations. For example, $\det{A}=1$ implies $B^{-1} = I\tr(B) - B$, by the Cayley-Hamilton theorem (or direct verification). Multiplying by $AB$ and taking traces gives
\begin{equation}
\tr(A) = \tr(AB)\tr(B) - \tr(AB^2).
\end{equation}
Thus $(\tr(A),\tr(B),\tr(AB))$ and $(\tr(AB^2),\tr(B),\tr(AB))$ are related by the Markoff move $m_1$.
To maintain the convention that the third matrix is the product of the first two, we note that $\tr(B) = \tr(B^{-1})$ and write the move as $A,B \mapsto AB^2, B^{-1}$. Writing an ``abelianized" vector that keeps track of the exponents on $A$ and $B$ but not the order of the product, we may write $m_1$ as a matrix
\[
[m_1] = \begin{bmatrix} 1 & 0 \\ 2 & -1 \end{bmatrix}
\]
exactly as in Proposition~\ref{prop:matrices}. Similar calculations changing the roles of $A, B$, and $AB$ give the matrices for the other moves $m_2$ and $m_3$. Likewise, the transpositions act by their corresponding matrices.
The key difference is that the action is no longer linear.

\section{Conclusion} \label{sec:conc}

We have investigated a family of 3-regular graphs defined from the solutions to (\ref{eqn:markoff}) modulo $p$ for each prime $p \geq 5$. It has already been conjectured that these graphs are connected \cite{Bar}, \cite{BGS}.
On the basis of the data summarized in Figure~\ref{fig:eigenvalue}, we further conjecture that these graphs are asymptotically Ramanujan for $p \equiv 3$ mod 4. That is, the second largest eigenvalue $\lambda_2(p)$ converges to $2\sqrt{2}$ in this case. For $p \equiv 1$ mod 4, we conjecture that $\lambda_2(p)$ converges to a limit strictly less than $3$ and larger than $2\sqrt{2}$, but we do not venture a guess as to its value. It seems that the limit is approximately 2.875... and that there are relatively few eigenvalues above $2\sqrt{2}$. Indeed, Figure~\ref{fig:outstanding_eigen} suggests that the number of exceptional eigenvalues is asymptotic to $C p$ for a constant $C > 0$. Gathering this data involved computing many eigenvalues instead of only $\lambda_2$, so we considered only an even smaller range of primes. Thus the value of $C$ may not be accurate, but we do conjecture that the exponent $p^1$ is correct, and in particular that these large eigenvalues comprise a vanishing proportion of the total of roughly $p^2$ eigenvalues. This means the bulk of the spectrum is supported on $[-2\sqrt{2}, 2\sqrt{2}]$ and we conjectured further that this distribution converges to the Kesten-McKay law. In the meantime, the Kesten-McKay law has now been verified theoretically \cite{DM}, and Figure~\ref{fig:histogram} already shows a good fit even for the small primes $p=83$ and $p=89$.

For the level surfaces with $k \neq 0$ in Equation~\ref{eqn:levelk}, connectedness is no longer guaranteed and the more basic question of how many components there are (that is, the multiplicity of $\lambda = 3$) replaces the finer spectral questions above. The most extreme case is when $p$ divides $9k-4$, and then the components can be understood in terms of a linear action and the Cayley cubic. In general, the component sizes are dictated by arithmetic relations between $k$ and $p$. The simplest example of this is that there is a component of size 6 whenever $k$ is a square modulo $p$.

\section*{Acknowledgments}
We thank Peter Sarnak for his advice, encouragement, and support over the course of our work. We thank Pedro Henrique Pontes for showing us Lemma~\ref{lem:pedro}, which provides a simpler way to count solutions than our original proof using Gauss sums. We thank ReMatch, a summer research program at Princeton University, for being a supportive and stimulating research community.
Lee was supported by the Bershadsky Family Summer Research Scholars Fund through ReMatch and the Office of Undergraduate Research at Princeton University. de Courcy-Ireland was supported by a PGS D grant from the Natural Sciences and Engineering Research Council of Canada.

\end{document}